\documentclass[12pt]{amsart}
\usepackage{amsmath,amssymb,amsthm,amsfonts,amsxtra}
\usepackage{color}

\newcommand{\changed}[1]{\textcolor{black}{#1}}
\newcommand{\pg}[1]{\textcolor{black}{#1}}

\renewcommand{\section}[1]{%
\bigbreak\vskip-\lastskip\bigbreak%
\stepcounter{section}\noindent%
{\bf\boldmath\large \S \arabic{section}.\; #1.}%
\nobreak\medskip\newline\nobreak}

\newcommand{\myauthor}[1]{\vspace{1truecm}\centerline{\large By {#1}}%
\bigskip}

\renewenvironment{thebibliography}[1]{%
\bigbreak\vskip-\lastskip\bigskip\bigbreak\centerline{\bf References}
\nobreak\bigskip\nobreak

\begin{enumerate}}
{\end{enumerate}}

\newtheorem{thm}{Theorem}[section]
\newtheorem{prop}[thm]{Proposition}
\newtheorem{lemma}[thm]{Lemma}

\theoremstyle{definition}
\newtheorem{dfn}[thm]{Definition}
\theoremstyle{remark}
\newtheorem{rem}[thm]{Remark}
\numberwithin{equation}{section}

\newcommand{\pf}{\textit{Proof. }}

\newcommand{\lbar}[1]{ \overline{#1} }
\newcommand{\norm}[1]{ \Vert #1 \Vert }
\newcommand{\set}[2]{ \left\{\,#1\,;\,#2\,\right\} }
\newcommand{\cpling}[2]{ \langle #1,#2 \rangle}
\newcommand{\transp}[1]{ {}^{\mathsf{t}} #1}

\newcommand{\nc}{\newcommand}

\nc{\ep}{\varepsilon}
\nc{\orbt}{\mathcal{O}}
\nc{\mychi}{\raise 1pt \hbox{$\chi$}}
\nc{\Pcone}{\mathcal{P}}
\nc{\Qcone}{\mathcal{Q}}
\nc{\Zsp}{\mathcal{Z}}
\nc{\Vsp}{\mathcal{V}}
\nc{\us}{\underline{s}}
\nc{\um}{\underline{m}}
\nc{\usigma}{\underline{\sigma}}
\nc{\uep}{\underline{\varepsilon}}
\nc{\uu}{\underline{u}}
\nc{\up}{\underline{p}}
\nc{\ud}{\underline{d}}
\nc{\meas}{\mathcal{M}}
\nc{\Rz}{\mathcal{R}}
\nc{\dirsum}{\sideset{\ }{^{\oplus}}\sum}
\nc{\iu}{i}
\nc{\integer}{\mathbb{Z}}
\nc{\real}{\mathbb{R}}
\nc{\complex}{\mathbb{C}}

\address{\noindent Piotr Graczyk 
\newline
    Laboratoire de Math\'ematiques LAREMA,  \newline
    Universit\'e d'Angers,
    2, boulevard Lavoisier,
    49045 Angers Cedex 01, France \vspace{5pt}}

\address{\noindent  Hideyuki Ishi 
\newline
  Graduate School of Mathematics, \newline
  Nagoya University,
  Furo-cho,
  464-8602 Nagoya, Japan
\vspace{10pt}}

\email{graczyk@univ-angers.fr, hideyuki@math.nagoya-u.ac.jp} 

\begin{document}
 
\setlength{\baselineskip}{18pt}
\centerline{\Large\bf Riesz measures and Wishart laws}
\bigskip
\centerline{\Large\bf associated to quadratic maps}
\myauthor{Piotr \textsc{Graczyk} and Hideyuki \textsc{Ishi}}
 
\bigskip

\footnotetext{
This research was partially supported by the grant  ANR-09-BLAN-0084-01.}

\begin{center}

\begin{minipage}{12truecm}\small
\textbf{Abstract.}
We introduce a natural definition of Riesz measures and Wishart laws
 associated to an $\Omega$-positive (virtual) quadratic map,
 where $\Omega \subset \real^n$ is a regular open convex cone.
We give a general formula for moments of the Wishart laws.
Moreover,
 if the quadratic map has an equivariance property
 under the action of a linear group 
 acting on the cone $\Omega$ transitively,
 then the associated Riesz measure and Wishart law are described explicitly
 by making use of 
 theory of relatively invariant distributions on homogeneous cones.
\end{minipage}
\end{center}

%
%
%
%
\section{Introduction}

Riesz measures and distributions on convex  cones  form  one of fundamental tools of harmonic analysis
and of the theory of the wave equation, cf. \cite{F-K}  in the case of  symmetric cones 
and \cite{G1}, \cite{I1} for homogeneous cones.   Moreover,
exponential families generated by  Riesz measures are composed of Wishart laws and  are
of great significance in   random matrix theory  and in statistics.\\

Wishart laws are probability distributions on 
symmetric or Hermitian matrices with very important applications
in   multivariate statistics. Their role in statistics
is due to two reasons:

-- they are probability distributions of the maximum likelihood estimator(MLE)
of the covariance matrix in a multivariate normal sample (\cite{Mu}, \cite{A-W}).

-- in Bayesian statistics, Wishart laws  form a Diaconis-Ylvisaker family (\cite{D-Y}) of prior distributions
for the covariance parameter in a covariance selection model (\cite{L-M}).

On the other hand, recent developments in random matrix theory  of chiral Gaussian   ensembles  containing
Wishart laws, are  intense  and motivated by applications in mathematical physics,
cf. \cite{K-T} and references therein. \\

These numerous modern applications of Wishart laws make it necessary to develop the theory of
 Wishart laws and Riesz measures on  more general cones than in the classical case of the symmetric cones
 of real symmetric or complex Hermitian matrices. 
 For example, in an $r$-dimensional Gaussian model \changed{$X$},  if the marginal variables \changed{$X_i$} and \changed{$X_j$}
 are known to be conditionally independent given all the other variables, the statistical
 analysis of the covariance matrix  of \changed{$X$}  must be done on the cone \changed{$\Pcone$}
 of positive definite symmetric matrices \changed{$Y$}  with \changed{$Y_{ij}=Y_{ji}=0$} and on its dual cone \changed{$\Qcone$}
 (\cite{L-M}).
 The cones $\Pcone$ and $\Qcone$ are usually no longer symmetric.  This led to  some important papers
 in recent statistical and probabilistic literature about  Wishart laws  on more general cones: homogeneous cones (\cite{A-W})
 or cones related to graphical models (\cite{L-M}).  In these papers, Wishart laws are introduced 
 via their density functions
 \changed{(see Section 3.8)}. \\
 
    In our paper we  construct and study  Riesz measures and  Wishart laws on regular convex cones
 via quadratic maps.  For a regular open convex cone
 $\Omega \subset \real^n$
 and an $\Omega$-positive quadratic map
 $q : \real^m \to \real^n$,
 the Riesz measure associated to $q$ is defined as the image of the Lebesgue measure $dx$ on $\real^m$  
 by $q$.   Wishart laws studied in this paper are obtained from    $\real^m$-valued normal random  vectors $X$ 
  as   the law  of $Y := q(X)$/2.
 This is a classical natural  approach to Riesz measures (\cite{F-K}) and   Wishart laws (\cite{Mu},\cite{F-K})
 and we propose to extend it to a much more general setting.\\
 
 In Section 2 of the paper we  explain the details of the quadratic construction of Riesz measures on regular convex cones
 and next we define the corresponding Wishart laws.  We compute their Laplace transforms, what is the starting
 point to get formulas for their expectation, covariance and higher  moments (Theorems \ref{thm:E_and_V} and
\ref {thm:moments}). More general Riesz and Wishart distributions associated to {\it virtual} quadratic maps are introduced in Section \changed{2.3}.
 Moments formulas are generalized (Theorem \ref{thm:gen_moments}).
 Group equivariance of Wishart laws is studied at the end of the section.
 
 Section 3 of the article is thoroughly devoted to the case when $\Omega$ is a homogeneous cone
 and  the quadratic map $q$  is homogeneous.  A crucial role in the analysis of these maps
 and of related Riesz measures and Wishart  laws  is played by a matrix realization of any
 homogeneous cone, coming from \pg{\cite{I5}} and explained in Section 3.2. It allows, among others,
 to define basic and standard quadratic maps in Sections 3.3 and 3.4. They play a role
 of generators for homogeneous quadratic \changed{maps} $q$ needed to construct all Riesz
 measures and Wishart laws on $\Omega$.
 Next we apply
 the results of \cite{I1} on Gindikin-Riesz distributions on $\Omega$ and \pg{on} the orbit decomposition of \changed{$\lbar{\Omega}$,
 the closure of $\Omega$}.
 We explain  the relation between   Riesz measures  related to homogeneous quadratic maps and the 
  Gindikin-Riesz distributions on $\Omega$ (Theorem \ref{thm:Riesz-Gindikin}).  
 \changed{In Section 3.7}, we prove the Bartlett decomposition for the Wishart laws on homogeneous cones
(Theorems \ref{thm:Bartlett1} and \changed{\ref{thm:Bartlett2}}).

Families of Wishart laws that we  construct and study  in Section 3 comprise 
Wishart distributions studied in papers \cite{A-W} and  \cite{L-M} (homogeneous case)
and are significantly bigger: we  describe  all singular Wishart laws and
many more absolutely continuous Wishart laws than in papers \cite{A-W} and \cite{L-M}.
\changed{For the symmetric cone case, our Wishart laws cover the ones studied in \cite{H-L} as well}. 
All the results of  Section \changed{2} apply to them, in particular the formulas for the moments.\\

Throughout the paper and from its very beginning all our concepts are illustrated on important
examples, including a non-homogeneous cone 
(Example 1 Section 2), symmetric cones of positive definite real symmetric matrices $\Pi_r$(Example 2 Section 2),
 4-dimensional Lorentz cone(Example 5 Section 3), and non-symmetric but homogeneous Vinberg
cone and its dual(Example 3 Section 2).\\[5mm]
{\bf Acknowledgement}. We thank Professors Gerard Letac and Yoshihiko Konno for discussions
on the topic of the article.

%
%
%
%
\section{Riesz measure and Wishart law on \changed{a} convex cone}
\noindent{\bf 2.1. ~Regular cones and quadratic maps.} \indent 
In this paper,
 an open convex cone
 $\Omega \subset \real^n$
 is always assumed to be \textit{regular},
 that is,
 $\lbar{\Omega} \cap (- \lbar{\Omega}) = \{0\}$,
 where $\lbar{\Omega}$ denotes the closure of $\Omega$.
Then
 the dual cone
 $\Omega^* 
  := \set{\eta \in (\real^n)^*}
         {\cpling{y}{\eta} > 0 \,\,
          (\forall y \in \lbar{\Omega} \setminus \{0\})}$
 is a regular open convex cone again 
 in the dual vector space
 $(\real^n)^*$,
 and we have $(\Omega^*)^* = \Omega$.
An $\real^n$-valued quadratic map
 $q : \real^m \to \real^n$ 
 is said to be \textit{$\Omega$-positive}
 if 
 (i) $q(x) \in \lbar{\Omega}$ for all $x\in \real^m$,
 and
 (ii) $q(x) = 0$ implies $x=0$.
These (i) and (ii) are restated in a single condition
 $q(x) \in \lbar{\Omega} \setminus \{0\} 
  \,\,\,(\forall x \in \real^m \setminus \{0\})$.
For the quadratic map $q$,
 we define \textit{the associated linear map}
 $\phi = \phi_{q}: (\real^n)^* \to \mathrm{Sym}(m, \real)$
 in such a way that
 $$
 \transp{x} \phi(\eta) x =  \cpling{q(x)}{\eta}
 \qquad (\eta \in (\real^n)^*,\,\,x \in \real^m).
 $$
Then the $\Omega$-positivity of $q$ is equivalent to
 the following property of $\phi$:
 \begin{equation} \label{eqn:phi_positive}
 \eta \in \Omega^* \Rightarrow
  \phi(\eta) \mbox{ is positive definite}. 
 \end{equation}

\noindent
{\bf Example 1.}\indent
Let $\Omega$ be the open convex cone in $\real^3$ defined by
\begin{align}
 \Omega 
 &:= \set{t_1 \begin{pmatrix} 0 \\ 0 \\ 1 \end{pmatrix}
         + t_2 \begin{pmatrix} 1 \\ 0 \\ 1 \end{pmatrix}
         + t_3 \begin{pmatrix} 1 \\ 1 \\ 1 \end{pmatrix}
         + t_4 \begin{pmatrix} 0 \\ 1 \\ 1 \end{pmatrix}
        }{t_1, \,t_2,\,t_3,\,t_4 >0} \nonumber\\
 &= \set{\begin{pmatrix} y_1 \\ y_2 \\ y_3 \end{pmatrix} \in \real^3}
        {y_1>0,\,\,y_2 >0,\,\,-y_1 + y_3 >0,\,\,-y_2 + y_3 >0}.
 \label{eqn:4dimcone}
\end{align}
If we identify $(\real^3)^*$ with $\real^3$ by
 $\cpling{y}{\eta} :=  y_1 \eta_1 + y_2 \eta_2 + y_3 \eta_3
 \,\,\,(y, \eta \in \real^3)$,
 we have
\begin{align*}
 \Omega^* 
 &= \set{t_1 \begin{pmatrix} 1 \\ 0 \\ 0 \end{pmatrix}
         + t_2 \begin{pmatrix} 0 \\ 1 \\ 0 \end{pmatrix}
         + t_3 \begin{pmatrix} -1 \\ 0 \\ 1 \end{pmatrix}
         + t_4 \begin{pmatrix} 0 \\ -1 \\ 1 \end{pmatrix}
        }{t_1, \,t_2,\,t_3,\,t_4 >0}\\
 &= \set{\begin{pmatrix} \eta_1 \\ \eta_2 \\ \eta_3 \end{pmatrix} \in \real^3}
        {\eta_3>0,\,\,\eta_1 + \eta_3>0,\,\,
         \eta_1 + \eta_2 + \eta_3 >0,\,\, \eta_2 + \eta_3 >0},
\end{align*}
 see \cite{I4}.
Let 
 $q : \real^4 \to \real^3$
 be the quadratic map given by
 $$
 q(x) := (x_1)^2 \begin{pmatrix} 0 \\ 0 \\ 1 \end{pmatrix}
         + (x_2)^2 \begin{pmatrix} 1 \\ 0 \\ 1 \end{pmatrix}
         + (x_3)^2 \begin{pmatrix} 1 \\ 1 \\ 1 \end{pmatrix}
         + (x_4)^2 \begin{pmatrix} 0 \\ 1 \\ 1 \end{pmatrix}
 \quad (x \in \real^4).
 $$
Clearly,
 this $q$ is $\Omega$-positive.
By a simple calculation, 
 we have
$$
 \phi(\eta) 
 = \begin{pmatrix} 
   \eta_3 & 0 & 0 & 0\\  0 & \eta_1 + \eta_3 & 0 & 0\\
   0 & 0 & \eta_1 + \eta_2 + \eta_3 & 0 \\ 0 & 0 & 0 & \eta_2 + \eta_3
  \end{pmatrix}
 \quad (\eta \in (\real^3)^*).   
$$
{\bf Example 2.}\indent
Let $\Pi_r$
 be the set of positive definite real symmetric matrices of
 size $r$.
Then $\Pi_r$ is a regular open convex cone 
 in the vector space $\mathrm{Sym}(r, \real)$
 of real symmetric matrices.
If we identify the space $\mathrm{Sym}(r, \real)$
 with its dual vector space by
 the inner product
 $\cpling{y}{\eta} := \mathrm{tr}\,(y \eta)
  \,\,\,(y,\,\eta \in \mathrm{Sym}(r, \real))$,
 then the dual cone
  $\Pi_r^*$ coincides with $\Pi_r$.
We define 
 $q_{r,s}: \mathrm{Mat}(r,s; \real)\to \mathrm{Sym}(r, \real)$
 by
 $$
 q_{r,s}(x) = x \, \transp{x}
 \quad (x \in \mathrm{Mat}(r,s; \real) ).
 $$
Then $q_{r,s}$ is $\Pi_r$-positive.
We denote the $(i,j)$ component of 
 $x \in \mathrm{Mat}(r,s; \real)$
 by $x_{r(j-1)+i}$,
 so that 
 $\mathrm{Mat}(r,s; \real)$ is identified with $\real^{rs}$.
Then we have 
 for $\eta \in \mathrm{Sym}(r, \real)$
 $$
 \phi(\eta) 
 = \begin{pmatrix} 
   \eta & & & \\ & \eta & & \\ & & \ddots & \\ & & & \eta
  \end{pmatrix}
 \in \mathrm{Sym}(rs, \real),
 $$
 where $\eta$ is put $s$ times.
In this case,
 the map 
 $\phi : \mathrm{Sym}(r, \real) \to \mathrm{Sym}(rs, \real)$
 is a Jordan algebra representation,
 and $q$ is exactly the quadratic map associated to the representation
 (\cite[Chapter IV, Section 4]{F-K}).\\
{\bf Example 3.}\indent
Let $\Zsp$ be a subspace of $\mathrm{Sym}(r, \real)$, 
 and put $\Pcone := \Zsp \cap \Pi_r$.
Then $\Pcone$ is a regular open convex cone in $\Zsp$.
Let $\Qcone \subset \Zsp^*$ be
 the dual cone of $\Pcone$.
We shall construct a $\Qcone$-positive
 quadratic map
 $q_{\Zsp} : \real^r \to \Zsp^*$ 
 whose associated linear map 
 $\phi_{\Zsp} : \Zsp \to \mathrm{Sym}(r, \real)$
 equals the inclusion map.
Let us define
 the surjective linear map 
 $\pi_{\Zsp^*} : \mathrm{Sym}(r, \real) \to \Zsp^*$ 
 by
 $$
 \cpling{y}{\pi_{\Zsp^*}(S)} := \mathrm{tr}\,y S
 \qquad (y \in \Zsp,\,\,S \in \mathrm{Sym}(r, \real)).
 $$
Then the quadratic map
 $q_{\Zsp} : \real^r \to \Zsp^*$ 
 is given by
 $q_{\Zsp}(x) := \pi_{\Zsp^*}(x \transp{x})
  \,\,\,(x \in \real^r)$.
In fact,
 for $x \in \real^r \setminus \{0\}$ and $y \in \Pcone$
 we have
 \begin{equation} \label{eqn:qZ}
 \cpling{y}{q_{\Zsp}(x)}
 = \mathrm{tr}\,(y x \transp{x})
 = \transp{x} y x >0
 \end{equation}
 because $y$ is positive definite.
Therefore we get
 $q_{\Zsp}(x) \in \lbar{\Qcone} \setminus \{0\}$,
 so that $q_{\Zsp}$ is $\Qcone$-positive.
Keeping 
 the natural isomorphism $(\Zsp^*)^* \simeq \Zsp$
 in mind,
 we see from (\ref{eqn:qZ}) that 
 $\phi_{\Zsp}(y) = y\,\,\,(y \in \Zsp)$.
Soon later,
 we shall consider the cases
\begin{equation} \label{eqn:red_cone}
\Zsp 
:= \set{ \begin{pmatrix} y_{11} & 0 & 0 \\ 0 & y_{22} & y_{32}\\
   0 & y_{32} & y_{33} \end{pmatrix} \in \mathrm{Sym(3, \real)} }
   {y_{11},\,y_{22},\,y_{32},\,y_{33} \in \real}
\end{equation}
and
\begin{equation} \label{eqn:dual-Vinberg}
\Zsp 
:= \set{ \begin{pmatrix} y_{11} & 0 & y_{31} \\ 0 & y_{22} & y_{32}\\
                         y_{31} & y_{32} & y_{33} \end{pmatrix} 
\in \mathrm{Sym(3, \real)} }
{y_{11},\,y_{21},\,y_{31},\,y_{22},\,y_{33} \in \real}
\end{equation}
 as concrete examples.
Actually, 
 in the latter case (\ref{eqn:dual-Vinberg}),
 the cones $\Qcone$ and $\Pcone$ are called 
 \textit{the Vinberg cone} and \textit{the dual Vinberg cone} respectively,
 which are the lowest dimensional non-symmetric homogeneous cones
 (\cite{V1}).
We shall see another realization of the Vinberg cone $\changed{\Qcone}$
 in (\ref{eqn:Vinberg_cone})
 \changed{and the last paragraph of Section 3.3.}

Let $I = \{i_1,\,i_2, \dots, i_{\changed{k}}\}$
 be a subset of $\{1, \dots, r \}$
 with $1 \le i_1 < i_2 < \dots < i_{\changed{k}} \le r$,
\changed{and define
\begin{equation} \label{eqn:def_of_RI}
 R^I := \set{x \in \real^r}{x_i = 0 \mbox{ if }i \notin I}.
\end{equation}
}
We denote by $q_{\Zsp}^I$ the restriction of $q_{\Zsp}$
 to the space $R^I \subset \real^r$.
Clearly $q_{\Zsp}^I : R^I \to \Zsp^*$ is $\Qcone$-positive.
The associated linear map
 $\phi_{q_{\Zsp}^I} : \Zsp \to \mathrm{Sym}(\changed{k}, \real)$
 gives a submatrix of elements $y \in \Zsp$,
 that is,
 $\phi_{q_{\Zsp}^I}(y) = (y^{}_{i_{\alpha}i_{\beta}})$,
 which we denote by $y^{}_I$.\\

\noindent{\bf 2.2. \pg{Riesz measures and }~Wishart laws associated to quadratic maps.}\indent
For a regular open convex cone
 $\Omega \subset \real^n$
 and an $\Omega$-positive quadratic map
 $q : \real^m \to \real^n$,
 let $\mu_q$ be the image of the Lebesgue measure $dx$ on $\real^m$  
 by $q$.
Namely, 
 the measure $\mu_q$ on $\real^n$ is defined
 in such a way that
 \begin{equation} \label{eqn:def_of_muq}
 \int_{\real^n} f(y) \,\mu_q(dy) = \int_{\real^m} f(q(x))\,dx
 \end{equation}
 for a measurable function $f$ on $\real^n$.
The $\Omega$-positivity of $q$ implies that
 the support of $\mu_q$
 is contained in the closure $\lbar{\Omega}$ of the cone $\Omega$.
By analogy to \cite[Proposition VII.2.4]{F-K},
 we call $\mu_q$ \textit{the Riesz measure} associated to $q$.


\begin{lemma} \label{lemma:Laplace_trans}
Let $\phi : (\real^n)^* \to \mathrm{Sym}(m, \real)$
 be the linear map associated to $q$.
Then, 
 for $\eta \in \Omega^*$,
 the Laplace transform
 $L_{\mu_q}(-\eta) 
  := \int_{\real^n} e^{-\cpling{\eta}{y}} \mu_q(dy)$ 
 of $\mu_q$
 equals
 $\pi^{m/2} (\det \phi(\eta))^{-1/2}$.
\end{lemma}

\pf
By definition,
 we have
 $L_{\mu_q}(-\eta) = \int_{\real^m} e^{-\transp{x}\phi(\eta)x} \,dx$.
Since $\phi(\eta)$ is positive definite, 
 the assertion follows from
 a formula of the Gaussian integral. 
\qed
$ $\\
  

\begin{dfn} \label{dfn:Wishart}
The members of the exponential family 
 $\{\gamma_{q,\theta}\}_{\theta \in -\Omega^*}$
 generated by $\mu_q$ 
 are called 
 \textit{the Wishart laws on $\Omega$ associated to $q$}.
Namely,
 \begin{equation} \label{eqn:def_of_Wishart}
 \gamma_{q,\theta}(dy) 
 := \frac{e^{\cpling{y}{\theta}}}{L_{\mu_q}(\theta)}
    \mu_q(dy)
 \qquad (y \in \real^n).
 \end{equation} 
\end{dfn}
$ $\\
 \begin{rem}\label{rem-random} 
By (\ref{eqn:def_of_muq}), (\ref{eqn:def_of_Wishart})
 and Lemma~\ref{lemma:Laplace_trans},
 we have for a measurable function $f$ on $\real^n$
$$
 \int_{\real^n} f(y) \gamma_{q, \theta}(dy)
 = 
 \pi^{-m/2} (\det \phi(-\theta))^{1/2}
 \int_{\real^m} f(q(x)) e^{-\transp{x}\phi(-\theta)x}\,dx.
$$
Putting $\Sigma := \phi(-\theta)^{-1}$
 and replacing the variable $x$ by $x/\sqrt{2}$,
 we rewrite the right-hand side as
 $$
 (2\pi)^{-m/2}(\det \Sigma)^{-1/2} \int_{\real^m} 
 f(q(x)/2) e^{-\transp{x}\Sigma^{-1} x/2}\,dx.
 $$
Therefore,
 if $X$ is an $\real^m$-valued random variable with 
 the normal law $N(0, \phi(-\theta)^{-1})$,
 then $\gamma_{q, \theta}$ is nothing else but 
 the law of $Y := q(X)$/2.
In particular,
 the classical Wishart law
 as defined in \cite[Definition 3.1.3]{Mu}
 coincides with our $\gamma_{q, \theta}$ in Example 2.
 \end{rem} 
%
%


\begin{prop} \label{prop:L_gamma}
Let $Y$ be an $\real^n$-valued random variable
 with the Wishart law $\gamma_{q, \theta}$.
Then the Laplace transform
 $L_{\gamma_{q,\theta}}(\eta) = E(e^{\cpling{Y}{\eta}})$
 of $\gamma_{q, \theta}$
 is given by
$$
L_{\gamma_{q,\theta}}(\eta)
 = \det (I_m + \phi(-\theta)^{-1} \phi(-\eta))^{-1/2}
$$
 for $\eta \in -\theta - \Omega^*$.
\end{prop}
\pf
By definition,
 we have
 $L_{\gamma_{q,\theta}}(\eta) = L_{\mu_q}(\theta)^{-1} L_{\mu_q}(\eta+\theta)$.
Thus the formula follows from 
 Lemma \ref{lemma:Laplace_trans}
 and the observation that
$(\det \phi(-\theta))^{-1} \det \phi(-\eta - \theta)
 = \det \bigl( \phi(-\theta)^{-1}(\phi(-\theta) + \phi(-\eta))\bigr)
 = \det (I_m + \phi(-\theta)^{-1} \phi(-\eta))$.
\qed
$ $\\
 
We shall consider the mean and the covariance 
 of the Wishart law $\gamma_{q, \theta}$.
First we fix a notation for the directional derivative of a function,
 and recall some basic formulas.
For a (vector-valued) function $f$ on a domain $U \subset \real^N$,
 a point $a \in U$ and a vector $v \in \real^N$,
 we denote by $D_v f(a)$ the directional derivative of $f$ at $a$
 given by
 $D_v f(a) := (\frac{d}{dh})_{h=0}f(a + h v).$
Now we consider some functions and their derivatives 
 on the domain $- \Pi_m$ in the vector space $\mathrm{Sym}(m, \real)$.
First we set
 $f(x) = \log \det (-x) \in \real$ 
 for $x \in -\Pi_m$.
Then we have
\begin{equation} \label{eqn:trace-inverse}
 D_v f(a) = -\mathrm{tr}\,(-a)^{-1}v
 \quad (a \in -\Pi_m,\,\,v \in \mathrm{Sym}(m, \real)). 
\end{equation}
Next we observe the case
 $f(x) = (-x)^{-1} \in \mathrm{Sym}(m, \real)$ for $x \in -\Pi_m$.
Then
 \begin{equation} \label{eqn:inverse-inverse}
 D_v f(a) = (-a)^{-1} v\, (-a)^{-1}
 \quad (a \in -\Pi_m,\,\,v \in \mathrm{Sym}(m, \real)). 
 \end{equation}


\begin{thm} \label{thm:E_and_V}
Let $Y$ be an $\real^n$-valued random variable
 with the Wishart law $\gamma_{q, \theta}$.\\
{\rm (i)}
For $\eta \in (\real^n)^*$,
 one has
 $$E(\cpling{Y}{\eta}) =
 \mathrm{tr}\,\phi(-\theta)^{-1}\phi(\eta)/2.$$
{\rm (ii)}
For $\eta, \eta' \in (\real^n)^*$,
 one has
$$
 E((\cpling{Y}{\eta} - M)(\cpling{Y}{\eta'} - M'))
 =
 \mathrm{tr}\,\phi(-\theta)^{-1}\phi(\eta)\phi(-\theta)^{-1}\phi(\eta')/2, 
 $$
 where
 $M := E(\cpling{Y}{\eta})$ 
 and
 $M' := E(\cpling{Y}{\eta'})$. 
\end{thm}
\pf
Since $\{\gamma_{q,\theta}\}_{\theta \in - \Omega^*}$
 is the exponential family generated by $\mu_q$,
 it is well known (\cite{L-S})
 that the mean $E(\cpling{Y}{\eta})$
 is given by the derivative
 $D_{\eta}\log L_{\mu}(\theta)$,
 while the covariance
 $E((\cpling{Y}{\eta} - M)(\cpling{Y}{\eta'} - M'))$
 equals
 $D_{\eta} D_{\eta'} \log L_{\mu}(\theta)$.
Thus the formulas follow from 
 Lemma~\ref{lemma:Laplace_trans},
 (\ref{eqn:trace-inverse}) and (\ref{eqn:inverse-inverse})
 because $\phi(\theta) \in -\Pi_m$.
\qed
$ $\\

Let us discuss higher moments of the Wishart law
 $\gamma_{q, \theta}$. 
The computation is reduced to 
 the derivations 
 of the Laplace transform
 $L_{\mu_q}(\theta) = \pi^{m/2} (\det \phi(-\theta))^{-1/2}$.
Namely we have for $\eta_1, \dots, \eta_N \in (\real^n)^*$
$$
 E(\cpling{Y}{\eta_1} \cpling{Y}{\eta_2} \dots \cpling{Y}{\eta_N})
 = \frac{D_{\eta_1} D_{\eta_2} \dots D_{\eta_N}L_{\mu_q}(\theta)}
        {L_{\mu_q}(\theta)}. 
$$ 
On the other hand,
 if $f(x) = \det(-x)^{-p}$ for $x \in - \Pi_m$, 
 where $p$ is a real constant,
 the derivative
 $D_{v_1} D_{v_2} \dots D_{v_N}f(a)\,\,\,
  (a \in -\Pi_m,\,\,
   v_1, \dots, v_N \in \mathrm{Sym}(m, \real))$
 is given in \cite[Lemma 5]{G-L-M1}.
We review the formula briefly.
For an element $\pi$ of the symmetric group $\mathfrak{S}_N$,
 we write $C(\pi)$ for the set of cycles of $\pi$.
For $y \in \Pi_m$, we denote by $r_{\pi}(y)$
 the multilinear form of $v_1, \dots, v_N$ given by
$$
r_{\pi}(y)(v_1, \dots, v_N) 
 := \prod_{c \in C(\pi)} 
    \mathrm{tr}\Bigl(\prod_{j \in c} y v_j \Bigr). 
$$ 
Then we have
$$
 D_{v_1} D_{v_2} \dots D_{v_N}f(a)
 = f(a) \cdot 
  \sum_{\pi \in \mathfrak{S}_N} 
   p^{\sharp C(\pi)} r_{\pi}((-a)^{-1})(v_1, \dots, v_N), 
$$ 
 which is deduced from (\ref{eqn:trace-inverse}) 
 and (\ref{eqn:inverse-inverse}) by induction
 (see \cite{G-L-M1}).
Making use of the formula,
 we obtain


\begin{thm} \label{thm:moments}
Let $Y$ be an $\real^n$-valued random variable
 with the Wishart law $\gamma_{q, \theta}$.
For $\eta_1,\, \eta_2, \dots, \eta_N \in (\real^n)^*$,
 one has
\begin{align*}
E(\cpling{Y}{\eta_1} \cpling{Y}{\eta_2} \dots \cpling{Y}{\eta_N})
 &= \sum_{\pi \in \mathfrak{S}_N} 
    \bigl(\frac{1}{2}\bigr)^{\sharp C(\pi)} 
    r_{\pi}(\phi(-\theta)^{-1})
           (\phi(\eta_1),\phi(\eta_2), \dots, \phi(\eta_N))\\
 &= \sum_{\pi \in \mathfrak{S}_N} 
    \bigl(\frac{1}{2}\bigr)^{\sharp C(\pi)}
    \prod_{c \in C(\pi)} 
    \mathrm{tr}\Bigl(\prod_{j \in c} \phi(-\theta)^{-1} \phi(\eta_j) \Bigr).
\end{align*}
\end{thm}


%
%
%
\noindent{\bf 2.3.~Wishart laws associated to virtual quadratic maps.}\indent
We shall consider virtual quadratic maps,
 that is a 'formal linear combination' of 
 quadratic maps, 
 and the associated Wishart laws.
First we introduce the notion of direct sum of quadratic maps.
Let $q_i : \real^{m_i} \to \real^n\,\,\,(i=1, \dots, s)$ 
 be $\Omega$-positive quadratic maps.
Then the direct sum
 $q = q_1 \oplus q_2 \oplus \dots \oplus q_s$
 is an $\real^n$-valued quadratic map on
 $\real^{m_1} \oplus \real^{m_2} \oplus \dots \oplus \real^{m_s}$
 given by
$$
 q(x) := q_1(x_1) + q_2(x_2) + \dots + q_s(x_s) 
 \quad \Bigl( x =\sum_{k=1}^s x_i,\,\,\,x_i \in \real^{m_i} \Bigr).
$$
It is easy to see that $q$ is also $\Omega$-positive.   
If $q_1 = q_2 = \dots = q_s$,
 then the direct sum $q$ is denoted by $q_1^{\oplus s}$.

The linear map 
 $\phi : (\real^n)^* \to \mathrm{Sym}(m, \real)
 \,\,\,(m := \sum_{i=1}^s m_i)$
 associated to the direct sum $q = \sum^{\oplus} q_i$
 is given by
 \begin{equation} \label{eqn:sum_phii}
 \phi(\eta) = 
  \begin{pmatrix} 
   \phi_1(\eta) & & & \\ & \phi_2(\eta) & & \\ & & \ddots & \\ & & &
  \phi_s(\eta)
  \end{pmatrix}
 \qquad (\eta \in (\real^n)^*).
 \end{equation}
Conversely,
 if a symmetric matrix $\phi(\eta)$ is expressed
  by $\phi_1(\eta), \dots, \phi_s(\eta)$ as above 
  for all $\eta \in (\real^n)^*$,
 then the corresponding quadratic map $q$ is 
 the direct sum of $q_1, \dots, q_s$.
In Example 1,
 the quadratic map $q:\real^4 \to \real^3$
 is the direct sum of 4 quadratic maps
 $q_i : \real \owns x \mapsto x^2 v_i \in \real^3\,\,(i=1,\dots,4)$,
 where
 $$
 v_1 := \begin{pmatrix} 0 \\ 0 \\ 1 \end{pmatrix}, \quad
 v_2 := \begin{pmatrix} 1 \\ 0 \\ 1 \end{pmatrix}, \quad
 v_3 := \begin{pmatrix} 1 \\ 1 \\ 1 \end{pmatrix}, \quad
 q_4 := \begin{pmatrix} 0 \\ 1 \\ 1 \end{pmatrix}.
 $$
In Example 2,
 we see that
 $q_{r,s}: \mathrm{Mat}(r,s; \real)\to \mathrm{Sym}(r, \real)$
 is naturally identified with $q_{r,1}^{\oplus s}$.
In Example 3 with $\Zsp$ given by (\ref{eqn:red_cone}), 
 we have $q_{\Zsp} = q_{\Zsp}^{ \{1\} } \oplus q_{\Zsp}^{ \{2,3\} }$,
 while we do not have such a decomposition for the case (\ref{eqn:dual-Vinberg}).

Let $q_i : \real^{m_i} \to \real^n\,\,(i=1,2)$ 
 be $\Omega$-positive quadratic maps,
 and $q$ the direct sum $q_1 \oplus q_2$.
Then it is easy to see that
 the measure $\mu_q$ equals the convolution 
 $\mu_{q_1} * \mu_{q_2}$.
Thus, 
 for $\theta \in - \Omega^*$
 we have
 $L_{\mu_q}(\theta) = L_{\mu_{q_1}}(\theta) L_{\mu_{q_2}}(\theta)$
 and
 $\gamma_{q, \theta} = \gamma_{q_1, \theta} * \gamma_{q_2, \theta}$.
In general,
 if we set
 $q = q_1^{\oplus s_1} \oplus q_2^{\oplus s_2} \oplus \dots
      \oplus q_t^{\oplus s_t}$ 
 for $\Omega$-positive quadratic maps
 $q_i : \real^{m_i} \to \real^n\,\,(i=1,2, \dots, t)$
 and positive integers $s_1,\,s_2, \dots, s_t$,
 then 
 we have
\begin{equation} \label{eqn:virtual_LaplaceTrans}
 L_{\mu_q}(\theta) = \prod_{i=1}^t L_{\mu_i}(\theta)^{s_i}
 \quad (\theta \in -\Omega^*).
\end{equation}

Now
 we remark that,
 even though $s_i$'s are not positive integers,
 there may exist a positive measure $\mu_q$ on $\lbar{\Omega}$
 for which the relation (\ref{eqn:virtual_LaplaceTrans}) holds.
In general,
 for real numbers $s_1, \dots, s_r$,
 we shall call a formal sum
 $q = q_1^{\oplus s_1} \oplus q_2^{\oplus s_2} \oplus \dots
      \oplus q_t^{\oplus s_t}$ 
 a \textit{virtual $\Omega$-positive quadratic map}.
The measure $\mu_q$ satisfying (\ref{eqn:virtual_LaplaceTrans})
 is called \textit{the Riesz measure associated to $q$.}
By the injectivity of the Laplace transform,
 the associated Riesz measure $\mu_q$ is unique if it exists.
In this case,
 the Wishart laws
 $\gamma_{q, \theta}\,\,(\theta \in -\Omega^*)$
 are defined again as members of the exponential family generated by
 $\mu_q$. 


\begin{prop} \label{prop:virtual_L_gamma}
Let $q_i : \real^{m_i} \to \real^n\,\,\,(i=1, \dots, t)$ be
 $\Omega$-positive quadratic maps.
Assume that there exists a measure $\mu_q$ 
 associated to the virtual quadratic map
 $q = q_1^{\oplus s_1} \oplus \dots
      \oplus q_t^{\oplus s_t}$
 for certain $s_1, \dots, s_t \in \real$.
Let $Y$ be an $\real^n$-valued random variable
 with the Wishart law $\gamma_{q, \theta}$.
Then
 the Laplace transform
 $L_{\gamma_{q,\theta}}(\eta) = E(e^{\cpling{Y}{\eta}})$
 of the law $\gamma_{q, \theta}$
 is given by
$$
L_{\gamma_{q,\theta}}(\eta)
 \changed{= \prod_{i=1}^t L_{\gamma_{q_i,\theta}}(\eta)^{s_i} }
 = \prod_{i=1}^t \det (I_{m_i} + \phi_i(-\theta)^{-1} \phi_i(-\eta))^{-s_i/2}
$$
 for $\eta \in -\theta - \Omega^*$.
\end{prop}

\changed{\pf 
We have
 $L_{\gamma_{q,\theta}}(\eta) = L_{\mu_q}(\theta)^{-1}L_{\mu_q}(\eta+\theta)$
 by definition,
 and the right-hand side equals
 $\prod_{i=1}^t \bigl( L_{\mu_i}(\theta)^{-1}L_{\mu_i}(\eta+\theta) \bigr)^{s_i}$
 by (\ref{eqn:virtual_LaplaceTrans}).
Since 
 $L_{\gamma_{q_i,\, \theta}}(\eta) = L_{\mu_i}(\theta)^{-1}L_{\mu_i}(\eta+\theta)$,
 we obtain the first equality.
The second equality follows from Proposition~\ref{prop:L_gamma}.
\qed}
$ $\\
 
Since we see immediately from (\ref{eqn:virtual_LaplaceTrans})
 that
 $$
 \log L_{\mu_q}(\theta)
 = \sum_{i=1}^t s_i \log L_{\mu_i}(\theta),
 $$
 the virtual versions of Theorem \ref{thm:E_and_V} is given as follows:


\begin{prop} \label{prop:virtual_E_and_V}
Under the same assumption of Proposition \ref{prop:virtual_L_gamma},
 one has\\
{\rm (i)}
 $E(\cpling{Y}{\eta}) =
 \sum_{i=1}^t s_i\, \mathrm{tr}\,\phi_i(-\theta)^{-1}\phi_i(\eta)/2$
 for $\eta \in (\real^n)^*$,\\
{\rm (ii)}
$E((\cpling{Y}{\eta} - M)(\cpling{Y}{\eta'} - M'))
 = \sum_{i=1}^t
 s_i\,\mathrm{tr}\,\phi_i(-\theta)^{-1}\phi_i(\eta)
                 \phi_i(-\theta)^{-1}\phi_i(\eta')/2$
 for $\eta, \eta' \in (\real^n)^*$,
 where
 $M := E(\cpling{Y}{\eta})$ 
 and
 $M' := E(\cpling{Y}{\eta'})$.  
\end{prop}


As for higher moments,
 we generalize the formula in Theorem \ref{thm:moments}
 as follows:

\begin{thm} \label{thm:gen_moments}
Under the same assumption of Proposition \ref{prop:virtual_L_gamma},
 one has
\begin{align}
{} & E(\cpling{Y}{\eta_1} \cpling{Y}{\eta_2} \dots \cpling{Y}{\eta_N})
 \nonumber \\
 &= \sum_{\pi \in \mathfrak{S}_N} 
    \bigl(\frac{1}{2}\bigr)^{\sharp C(\pi)}
    \prod_{c \in C(\pi)} 
    \left\{ \sum_{i=1}^t s_i\,
    \mathrm{tr}\Bigl(\prod_{j \in c} \phi_i(-\theta)^{-1} \phi_i(\eta_j) \Bigr) \right\}.
 \label{eqn:gen_moments}
\end{align}
for $\eta_1,\, \eta_2, \dots, \eta_{N} \in (\real^n)^*$.
\end{thm}

Theorem~\ref{thm:gen_moments} easily follows from
 Theorem \ref{thm:moments}
 when $s_1, \dots, s_t$ are positive integers,
 that is, $q$ is a true quadratic map.
Indeed,
 we can consider the associated linear map $\phi = \phi_q$ in this case,
 and we have
 $$
 \mathrm{tr}\Bigl(\prod_{j \in c} \phi(-\theta)^{-1} \phi(\eta_j) \Bigr) 
 =
 \sum_{i=1}^t s_i\,
    \mathrm{tr}\Bigl(\prod_{j \in c} \phi_i(-\theta)^{-1} \phi_i(\eta_j) \Bigr)
 $$
 by virtue of (\ref{eqn:sum_phii}).
To prove (\ref{eqn:gen_moments}) for general case,
 it is enough to verify that  
 the quantity $E(\cpling{Y}{\eta_1} \cpling{Y}{\eta_2} \dots \cpling{Y}{\eta_N})$
 is a polynomial of $s_1, \dots, s_t$.
For this purpose, 
 we make some calculations involving 
 the semi-invariants or the cummulants
 (cf. \cite{L-S}).

For $\eta \in (\real^n)^*,\,\,\theta \in - \Omega^*$
 and a positive integer $k$,
 we denote by $S_k(\theta; \eta)$ 
 the $k$-th derivative of $\log L_{\mu_q}(\theta)$
 in the direction $\eta$:
$$
 S_k(\theta; \eta) = (D_{\eta})^k \log L_{\mu_q}(\theta)
 = \Bigl( \frac{d^k}{dh^k} \Bigr)_{h=0} \log L_{\mu_q}(\theta + h \eta).
$$
Similarly to Proposition~\ref{prop:virtual_E_and_V},
 we get
 \begin{equation} \label{eqn:Sk}
 S_k(\theta; \eta) = \frac{(k-1)!}{2}\sum_{i=1}^t
 s_i\,\mathrm{tr}\,\Bigl(\phi_i(-\theta)^{-1}\phi_i(\eta) \Bigr)^k.
 \end{equation}
On the other hand,
 since the function $L_{\mu_q}(\theta)$ is analytic,
 we have
 $$
 \log L_{\mu_q}(\theta + h \eta) 
 = \log L_{\mu_q}(\theta)+ \sum_{k=1}^{\infty} \frac{h^k}{k!} S_k(\theta;\eta).
 $$
Taking the exponential,
 we have
 \begin{equation} \label{eqn:Lmuq1}
 \begin{aligned}
 L_{\mu_q}(\theta + h \eta)
 &= L_{\mu_q}(\theta) \exp\Bigl(\sum_{k=1}^{\infty} \frac{h^k}{k!} S_k(\theta;\eta) \Bigr)\\
 &= L_{\mu_q}(\theta)\sum_{\ell=0}^{\infty} \frac{1}{\ell !}
   \Bigl(\sum_{k=1}^{\infty} \frac{h^k}{k!} S_k(\theta;\eta)\Bigl)^{\ell}.
 \end{aligned}
 \end{equation}
Let $Y$ be an $\real^n$-valued random variable
 with the Wishart distribution $\gamma_{q, \theta}$.
Then we have
 $$ 
 E(\cpling{Y}{\eta}^N) = 
  \frac{(D_{\eta})^N L_{\mu_q}(\theta)}
        {L_{\mu_q}(\theta)},
 $$
 so that
 \begin{equation} \label{eqn:Lmuq2}
 L_{\mu_q}(\theta + h \eta )
 = L_{\mu_q}(\theta) 
   \sum_{N=0}^{\infty} \frac{h^N}{N!} E(\cpling{Y}{\eta}^N).
 \end{equation}
Comparing the coefficients of $h^N$ 
 in (\ref{eqn:Lmuq1}) and (\ref{eqn:Lmuq2}),
 we obtain
\begin{align*}
 E(\cpling{Y}{\eta}^N) 
  = \sum_{\ell=1}^N \frac{1}{\ell!}\sum_{k_1 + k_2 + \dots + k_{\ell} = N} 
  \frac{N!}{k_1 ! k_2 ! \cdots k_{\ell}!}
  S_{k_1}(\theta; \eta) S_{k_2}(\theta; \eta)
  \dots S_{k_\ell}(\theta; \eta).
 \end{align*}
Substituting (\ref{eqn:Sk}) to the equality above,
 we obtain a formula for the higher moments.


\begin{prop} \label{prop:general_moment}
One has
 \begin{align} \label{eqn:eta-N}
 E(\cpling{Y}{\eta}^N) 
  = \sum_{\ell=1}^N \frac{1}{\ell!}\sum_{k_1 + k_2 + \dots + k_{\ell} = N} 
  \frac{N!}{k_1  k_2  \cdots k_{\ell}}
  \prod_{j=1}^{\ell} \Bigl(
 \sum_{i=1}^t
 \frac{s_i}{2}\,\mathrm{tr}\,\bigl(\phi_i(-\theta)^{-1}\phi_i(\eta) \bigr)^{k_j}
 \Bigr).
 \end{align}
\end{prop}


The moments
 $E(\cpling{Y}{\eta_1} \cpling{Y}{\eta_2} \cdots \cpling{Y}{\eta_N})
 \,\,\,(\eta_1, \dots, \eta_N \in (\real^n)^*)$
 can be obtained now by polarization,
 so that it is a polynomial of $s_1, \dots, s_t$.
Hence the proof of Theorem~\ref{thm:gen_moments} is completed.
We note that (\ref{eqn:gen_moments}) with $\eta_1 = \cdots = \eta_N = \eta$
 becomes (\ref{eqn:eta-N}),
 which means that 
 Proposition \ref{prop:general_moment} together with the polarization process
 also yields Theorem~\ref{thm:gen_moments}.\\

\noindent{\bf \changed{2.4}.~Group equivariance of the Wishart laws.} \indent
Let $G(\Omega)$ be the linear automorphism group 
 $\set{g \in GL(n, \real)}{g \Omega = \Omega}$ of $\Omega$.
For an $\Omega$-positive quadratic map $q : \real^m \to \real^n$
 and $g \in G(\Omega)$,
 the quadratic map $g \circ q : \real^m \to \real^n$ is again $\Omega$-positive.
It is easy to see that
 the Riesz measure $\mu_{g\circ q}$ is the image of $\mu_q$ by $g$,
 that is,
\begin{equation} \label{eqn:mugq}
 \mu_{g \circ q}(A) = \mu_q(g^{-1}A)
\end{equation}
 for a measurable set $A \subset \real^n$.
Let us discuss
 the Wishart laws $\gamma_{g \circ q, \theta}$ for $\theta \in - \Omega^*$.
For $\eta \in (\real^n)^*$,
 we denote by $g^*\eta$ the linear form $\eta \circ g \in (\real^n)^*$.
If $\eta \in \Omega^*$, 
 then
 $g^* \eta \in \Omega^*$
 because $\cpling{y}{g^*\eta} = \cpling{g y}{\eta} >0$
 for $y \in \lbar{\Omega} \setminus \{0\}$.
We observe
\begin{equation} \label{eqn:Lmugq}
 L_{\mu_q}(g^* \theta) = \int_{\real^m} e^{\cpling{q(x)}{g^*\theta}} dx =
 \int_{\real^m} e^{\cpling{g \circ q(x)}{\theta}} dx = L_{\mu_{g \circ q}}(\theta). 
\end{equation}
Therefore,
 denoting by $1_A$ the characteristic function of
 a measurable set $A \subset \real^n$,
 we have 
\begin{equation} \label{eqn:gamma-g}
\begin{aligned}
 \gamma_{g \circ q, \theta}(A)
 &= \frac{1}{L_{\mu_{g \circ q}}(\theta)}
 \int_{\real^m} 1_A(g \circ q(x)) e^{\cpling{g \circ q(x)}{\theta}} dx\\
 &=\frac{1}{L_{\mu_{q}}(g^*\theta)}
 \int_{\real^m} 1_{g^{-1}A}(q(x)) e^{\cpling{q(x)}{g^*\theta}} dx 
 = \gamma_{q, g^* \theta}(g^{-1}A).
\end{aligned}
\end{equation}
\changed{We restate (\ref{eqn:gamma-g}) as follows}.


\begin{lemma} \label{lemma:trans_gamma}
\changed{Let $g$ be an element of $G(\Omega)$.}
If a random variable $Y$ obeys the Wishart law $\gamma_{q, \theta}$, 
 the law of $gY$ is $\gamma_{g \circ q,\, (g^{-1})^*\theta}$.
\end{lemma}


\changed{
Let $q_i : \real^{m_i} \to \real^n\,\,\,(i=1, \dots, t)$ be
 $\Omega$-positive quadratic maps,
 and $q$ the virtual quadratic map
 $q_1^{\oplus s_1} \oplus \dots
      \oplus q_t^{\oplus s_t}$
 with $s_1, \dots, s_t \in \real$.
Then 
 we define $g \circ q$ to be the virtual quadratic map
 $(g \circ q_1)^{\oplus s_1} \oplus \dots
      \oplus (g \circ q_t)^{\oplus s_t}$.}

\changed{
\begin{prop} \label{prop:virtual_trans}
If the Riesz measure $\mu_q$ exists, 
 then 
 the Riesz measure $\mu_{g \circ q}$ exists and 
 equals the image of $\mu_q$ by $g$.
Moreover $\gamma_{g \circ q,\,(g^{-1})^*\theta}$
 is the image of $\gamma_{q, \theta}$ by $g$.
\end{prop}
\pf
Let $\mu'$ be the image of $\mu_q$ by $g$.
For $\theta \in - \Omega^*$,
 we have
 \begin{equation} \label{eqn:Lmuprime}
 L_{\mu'}(\theta) = \int_{\real^n} e^{\cpling{y}{\theta}} \mu'(dy) 
 = \int_{\real^n} e^{\cpling{g y}{\theta}} \mu_q(dy) 
 = \int_{\real^n} e^{\cpling{y}{g^* \theta}} \mu_q(dy) 
 = L_{\mu_q}(g^*\theta). 
 \end{equation}
By (\ref{eqn:virtual_LaplaceTrans}) and (\ref{eqn:Lmugq}),
 the last term equals
 $\prod_{i=1}^t L_{\mu_i}(g^*\theta)^{s_i} = \prod_{i=1}^t L_{\mu_{g \circ q_i}}(\theta)^{s_i}$.
Thus we get
 $L_{\mu'}(\theta)=\prod_{i=1}^t L_{\mu_{g \circ q_i}}(\theta)^{s_i}$
 which means $\mu' = \mu_{g \circ q}$ by (\ref{eqn:virtual_LaplaceTrans}).}

\changed{
Let $\gamma'$ be the image of $\gamma_{q, \theta}$ by $g$.
Similarly to (\ref{eqn:Lmuprime}),
 we have 
 $L_{\gamma'}(\eta) = L_{\gamma_{q, \theta}}(g^*\eta)$
 for $\eta \in - (g^{-1})^*\theta - \Omega^*$,
 while
 $L_{\gamma_{g \circ q_i, (g^{-1})^*\theta}}(\eta) = L_{\gamma_{q_i, \theta}}(g^*\eta)$
 by Lemma \ref{lemma:trans_gamma}.
On the other hand,
 we see from Proposition \ref{prop:virtual_L_gamma} 
 that
 $$L_{\gamma_{q, \theta}}(g^*\eta) = \prod_{i=1}^t L_{\gamma_{q_i,\theta}}(g^*\eta)^{s_i}
   = \prod_{i=1}^t L_{\gamma_{g \circ q_i, (g^{-1})^*\theta}}(\eta)^{s_i} = L_{\gamma_{g \circ q, (g^{-1})^*\theta}}(\eta).$$
Thus we get
 $L_{\gamma'}(\eta) = L_{\gamma_{g \circ q, (g^{-1})^*\theta}}(\eta)$,
 so that $\gamma' = \gamma_{g \circ q, (g^{-1})^*\theta}$ by the injectivity of the Laplace transform.
\qed
}
$ $\\
 
Let $\mathrm{Aut}(\Omega, q)$
 be the set of pairs
 $(g_1, g_2) \in G(\Omega) \times GL(m,\real)$
 for which
 $g_1 \circ q = q \circ g_2$. 
Then
 $\mathrm{Aut}(\Omega, q)$ forms
 a Lie subgroup of 
 $GL(n, \real) \times GL(m, \real)$,
 and we have a group homomorphism
$$
 pr_1 : \mathrm{Aut}(\Omega, q) \owns (g_1, g_2) \mapsto g_1 \in G(\Omega).
$$
The condition $(g_1,g_2) \in \mathrm{Aut}(\Omega, q)$
 is also equivalent to
 \begin{equation} \label{eqn:dualg}
 \phi_q(g_1^* \eta) = \transp{g_2} \phi_q(\eta) g_2
 \quad (\eta \in (\real^n)^*).
 \end{equation}
Then we obtain
\begin{equation} \label{eqn:RelInv}
 \det \phi_q(g_1^*\eta) = C \det \phi_q(\eta) \quad (\eta \in (\real^n)^*)
\end{equation}
 with $C = (\det g_2)^2$, 
 which means that 
 $\det \phi_q(\eta)$ is a relatively invariant polynomial on $(\real^n)^*$
 under the contragredient action of $pr_1(\mathrm{Aut}(\Omega, q))$.
The following proposition describes
 a transformation rule of the family of the Wishart laws $\{\gamma_{q,\theta}\}_{\theta \in -\Omega^*}$
 under the group $pr_1(\mathrm{Aut}(\Omega, q))$.
 

\begin{prop} \label{prop:trans}
For a measurable set $A \subset \real^n$
 and $(g_1, g_2) \in \mathrm{Aut}(\Omega, q)$, 
 one has\\
{\rm (i)}
 $\mu_q(g_1^{-1}A) = \mu_{g_1 \circ q}(A) = |\det g_2|^{-1} \mu_q(A)$,\\
{\rm (ii)} $\gamma_{q, g_1^*\theta}(A) = \gamma_{q, \theta}(g_1 A)$.
\end{prop}

\pf
{\rm (i)}
Because of (\ref{eqn:mugq}),
 we only have to show the second equality.
By definition, we have
$$
 \mu_{g_1 \circ q}(A) 
 = \int_{\real^m} 1_A(g_1 \circ q(x))\,dx
 = \int_{\real^m} 1_A(q \circ g_2(x))\,dx
$$
Putting $x' = g_2 x$,
 the last term equals
$$|\det g_2|^{-1} \int_{\real^m} 1_A(q(x'))\,dx'
 = |\det g_2|^{-1} \mu_q(A),$$  
whence (i) follows.\\
{\rm (ii)}
By (\ref{eqn:gamma-g}),
 we get for $y \in \real^n$
$$
\gamma_{q, g_1^*\theta}(dy) 
 = \gamma_{g_1 \circ q, \theta}(g_1 dy)  
 = \frac{e^{\cpling{g_1 y}{\theta}}}{L_{\mu_{g_1 \circ q}}(\theta)}\mu_{g_1 \circ q}(g_1dy),
$$
Since $\mu_{g_1 \circ q} = |\det g_2|^{-1} \mu_q$ by (i),
 the last term equals
$$
 \frac{e^{\cpling{g_1 y}{\theta}}}{|\det g_2|^{-1} L_{\mu_q}(\theta)}|\det g_2|^{-1} \mu_{q}(g_1 dy)\\
 = \gamma_{q, \theta}(g_1 dy).  
$$
Hence (ii) is verified.
\qed
$ $\\

%
%
%
%
\section{Homogeneous Case}
\noindent{\bf 3.1. ~Homogeneous quadratic map.}\indent
An $\Omega$-positive map $q : \real^m \to \real^n$
 is said to be \textit{homogeneous}
 if, for any $y,\,y' \in \Omega$,
 there exists $(g_1, g_2) \in \mathrm{Aut}(\Omega, q)$ 
 for which $g_1 y = y'$.
In other words,
 $q$ is homogeneous if $pr_1(\mathrm{Aut}(q,\Omega))$ acts on $\Omega$ transitively.
In this case,
 $\Omega$ is clearly a homogeneous cone,
 that is,
 a linear group on $\real^n$ acts on the cone $\Omega$ transitively.
\changed{
Then the dual cone $\Omega^* \subset (\real^n)^*$ is also a homogeneous cone
 on which the group
 $pr_1(\mathrm{Aut}(q,\Omega))$ acts transitively by the contragredient action
 (\cite{V1}). 
We see from (\ref{eqn:phi_positive}) and (\ref{eqn:dualg}) that
 the quadratic map $q$ is homogeneous 
 if and only if the associated linear map 
 $\phi_q: (\real^n)^* \to \mathrm{Sym}(m,\real)$
 is a representation of the dual cone $\Omega^*$
 in the sense of Rothaus \cite{Ro}
 (see also \cite{I7}).
}

A typical example of a homogeneous cone is
 $\Pi_r \subset \mathrm{Sym}(r, \real)$.
For $A \in GL(r,\real)$, 
 we denote by $\rho(A)$
 the linear map on $\mathrm{Sym}(r,\real)$ 
 defined by 
 $\rho(A)y := A y\, \transp{A}
  \,\,\,(y \in \mathrm{Sym}(r,\real))$.
Then the group $\rho(GL(r,\real))$
 acts on $\Pi_r$ transitively.
\changed{Moreover},
 the linear automorphism group $G(\Pi_r)$ equals $\rho(GL(\changed{r},\real))$.
We see that
 the quadratic map
 $q_{r,s}: \mathrm{Mat}(r,s; \real)\to \mathrm{Sym}(r, \real)$
 in Example 2 is homogeneous.
Indeed,
 we have a surjective homomorphism
 $GL(r,\real) \times O(s) \owns (A,B) \mapsto 
  (\rho(A),\,\tau_{r,s}(A,B)) \in \mathrm{Aut}(\Pi_r,\,q_{r,s})$,
 where $\tau_{r,s}(A,B)$ is a linear map on 
 $\mathrm{Mat}(r,s; \real)$
 given by
 $\tau_{r,s}(A,B)x := A x B^{-1}
  \,\,\,(x \in \mathrm{Mat}(r,s; \real))$.
\changed{
Thus we have 
 $pr_1(\mathrm{Aut}(q_{r,s},\Pi_r)) = \rho(GL(r,\real))$,
 which acts on $\Pi_r$ transitively.} 

For a subset $I \subset \{1, \dots, r\}$,
 we denote by $q^I$ the restriction of
 $q_{r,1} : \real^r \to \mathrm{Sym}(r,\real)$
 to the space $R^I \subset \real^r$
\pg{defined }  \changed{in (\ref{eqn:def_of_RI})}.
The map $q^I$ coincides with $q^I_{\Zsp}$ in Example 3
 with $\Zsp = \mathrm{Sym}(r,\real)$.
Let us observe that
 $q^I$ is homogeneous in general.
Let $P^I$ be the linear group
 consisting of $A \in GL(r,\real)$
 for which $A R^I = R^I$.
For example, if $r=3$, we have
\begin{align}
 P^{\{1\}} 
 &= \set{A = \begin{pmatrix} a_{11} & a_{12} & a_{13} \\
 0 & a_{22} & a_{23} \\ 0 & a_{32} & a_{33}\end{pmatrix}} 
 {A \in GL(3,\real)},
 \label{eqn:P1}\\
P^{\{2,3\}} 
 &= \set{A = \begin{pmatrix} a_{11} & 0 & 0 \\
 a_{21} & a_{22} & a_{23} \\ a_{31} & a_{32} & a_{33}\end{pmatrix}}
 {A \in GL(3,\real)}.
 \label{eqn:P23}
\end{align}

Since we have a homomorphism
 $P^I \owns A \mapsto (\rho(A),\,A) \in 
  \mathrm{Aut}(\Pi_r,\,q^I)$,
 it is enough to show that $\rho(P^I)$ acts on $\Pi_r$ transitively.
Put $k := \sharp I$ and
 take a permutation matrix
 $w_0 \in \mathfrak{S}_r \subset GL(r, \real)$
 sending
 $R^{\{r-k+1, \dots, r\}}$ onto $R^I$.
Then we have
$P^I = w_0 P^{\{r-k+1, \dots, r\}} w_0^{-1}$,
 and
$$
 P^{\{r-k+1, \dots, r\}} 
 = \set{\begin{pmatrix} A_1 & 0 \\ A_2 & A_3 \end{pmatrix}}
  {\begin{aligned} A_1 &\in GL(k, \real),\quad A_2 \in \mathrm{Mat}(r-k+1,k;\real)\\
   A_3 &\in GL(r-k+1,\real) \end{aligned}}. 
$$
Since $P^{\{r-k+1, \dots, r\}}$ contains the group of lower triangular matrices, 
 $\rho(P^{\{r-k+1, \dots, r\}})$ acts on $\Pi_r$ transitively.
Therefore $\rho(P^I) = \rho(w_0) \rho(P^{\{r-k+1, \dots, r\}}) \rho(w_0)^{-1}$ also acts on $\Pi_r$ transitively,
 so that $q^I$ is homogeneous.

Coming back to the examples (\ref{eqn:P1}) and (\ref{eqn:P23}),
 we note that $q = q^{\{1\}}  \oplus q^{\{2,3\}}$ is not homogeneous
 as $\Pi_3$-positive quadratic map,
 while both $q^{\{1\}}$ and $q^{\{2,3\}}$ are.
Indeed,
 the image of $q$ generates the space
 $\Zsp \subset \mathrm{Sym}(3,\real)$
 in (\ref{eqn:red_cone}).
Thus,
 if $(g_1, g_2) \in \mathrm{Aut}(\Pi_3,q)$,
 then $g_1$ must preserve both $\Zsp$ and $\Pi_3$.
Let us take $y \in \Pi_3 \setminus \Zsp$.
Then $g_1$ \changed{does not} send $I_3 \in \Pi_3$ to $y$ because $I_3 \in \Zsp$.
Thus the action of $pr_1(\mathrm{Aut}(q,\Pi_3))$ on $\Pi_3$ is not transitive. 

On the other hand,
 if we regard $q$ as a map from $R^{\{1\}} \oplus R^{\{2,3\}}$ to $\Zsp$,
 then $q$ is a homogeneous $\Pcone$-positive quadratic map,
 where $\Pcone := \Zsp \cap \Pi_r$.
In fact,
 since $(\rho(A), A) \in \mathrm{Aut}(q^{\{1\}},\Pcone) \cap \mathrm{Aut}(q^{\{2,3\}},\Pcone)$
 for $A \in P^{\{1\}} \cap P^{\{2,3\}}$,
 we have $\rho(A) \in pr_1(\mathrm{Aut}(q,\Pcone))$.
Therefore $pr_1(\mathrm{Aut}(q,\Pcone))$ contains
 a group $\rho(P^{\{1\}} \cap P^{\{2,3\}})$
 which acts on $\Pcone$ transitively.
  
In Example 1,
 the quadratic map $q : \real^4 \to \real^3$ is not homogeneous 
 because $\Omega \subset \real^4$ in (\ref{eqn:4dimcone})
 is not a homogeneous cone (\cite{I4}).\\

\noindent{\bf 3.2. ~Matrix realization of homogeneous cones.}\indent
In this section,
 we shall discuss a homogeneous cone
 realized as $\Pcone_{\Vsp} = \Zsp_\Vsp \cap \Pi_N$
 with 
 $\Zsp_{\Vsp} \subset \mathrm{Sym}(N, \real)$ 
 constructed from an appropriate system
 $\mathcal{V} = \{\mathcal{V}_{lk}\}$ of vector spaces
 in a specific way explained below,
 following \cite[section 3.1]{I5}.
The investigation of such cones
 is fundamental because all homogeneous cones are linearly equivalent 
 to some $\Pcone_{\Vsp}$ due to \cite[Theorem D]{I5}.

Let us take a partition
 $N = n_1+ n_2 + \dots+ n_r$ of a positive integer $N$,
 and consider a system of vector spaces
 $\Vsp_{lk} \subset \mathrm{Mat}(n_l,n_k;\real) \,\,\,(1 \le k < l \le r)$
 satisfying the following three conditions:\\
(V1)
 $A \in \Vsp_{lk},\,\, B \in \Vsp_{kj} \Rightarrow AB \in \Vsp_{lj}
  \quad (1 \le j < k < l \le r),$\\
(V2)
 $A \in \Vsp_{lj},\,\, B \in \Vsp_{kj} \Rightarrow A\,\transp{\! B} \in \Vsp_{lk}
  \quad (1 \le j < k < l \le r),$\\
(V3)
 $A \in \Vsp_{lk} \Rightarrow A\, \transp{\! A} \in \real I_{n_l}
  \quad (1 \le k < l \le r).$\\
Let $\Zsp_{\Vsp}$ be the subspace of 
 $\mathrm{Sym}(N, \real)$
 defined by
 $$
 \Zsp_{\Vsp}
 := \set{
 y = \begin{pmatrix}
 Y_{11} & \transp{Y_{21}} & \cdots & \transp{Y_{r1}} \\
 Y_{21} & Y_{22} & & \transp{Y_{r2}} \\
 \vdots & & \ddots & \\
 Y_{r1} & Y_{r2} & \cdots & Y_{rr} \end{pmatrix} 
 }{\begin{aligned}
 Y_{kk} &= y_{kk} I_{n_k},\,\,y_{kk} \in \real \,\,\,(k=1, \dots, r)\\
 Y_{lk} &\in \Vsp_{lk} \,\,\,(1 \le k < l \le r)
 \end{aligned} }.
$$ 
We set
 $\Pcone_{\Vsp} := \Zsp_\Vsp \cap \Pi_N$.
Then
 $\Pcone_{\Vsp}$
 is a regular open convex cone in the vector space $\Zsp_\Vsp$.
Let $H_N$ be the group of real lower triangular matrices with positive 
 diagonals,
 and $H_\Vsp$ a Lie subgroup of $H_N$
 defined by
 $$ 
H_\Vsp := \set{
 T = \begin{pmatrix}
 T_{11} & & & \\
 T_{21} & T_{22} & & \\
 \vdots & & \ddots & \\
 T_{r1} & T_{r2} & \cdots & T_{rr} \end{pmatrix} 
 }{\begin{aligned}
 T_{kk} &= t_{kk} I_{n_k},\,\,t_{kk} >0 \,\,\,(k=1, \dots, r)\\
 T_{lk} &\in \Vsp_{lk} \,\,\,(1 \le k < l \le r)
 \end{aligned}
 }.
$$ 
If $T \in H_\Vsp$ and $y \in \Zsp_\Vsp$, 
 then $\rho(T)y = T y \,\transp{T}$ \changed{belongs to} $\Zsp_\Vsp$
 thanks to (V1)--(V3).
Moreover $\rho(H_\Vsp)$ acts on the cone $\Pcone_\Vsp \subset \Zsp_\Vsp$
 simply transitively (cf. \cite[Proposition 2.1]{I5}).

Keeping (V3) in mind, 
 we define an inner product on the vector space $\Vsp_{lk}\,\,\,(1 \le k < l \le r)$ 
 by
 the equality 
\begin{equation} \label{eqn:innerprod_Vlk}
 A \,\transp{\! A} = (A|A) I_{n_l} = \norm{A}^2 I_{n_l} \quad (A \in \Vsp_{lk}).
\end{equation}
For $y, y' \in \Zsp_\Vsp$, 
 we set
\begin{equation}
 \cpling{y}{y'} := \sum_{k=1}^r y_{kk} y'_{kk} + 2 \sum_{1 \le k < l \le r} (Y_{lk}|Y'_{lk}), 
\end{equation}
 where $y_{kk}$ and $Y_{lk}$ (respectively $y'_{kk}$ and $Y'_{lk}$) denote the components of $y$ (respectively $y'$).
Note that the \changed{inner} product is not equal to $\mathrm{tr}\,y y'$ 
 unless $n_1 = \cdots = n_r = 1$.
By this coupling,
 we identify the dual space $\Zsp^*_\Vsp$ with $\Zsp_{\Vsp}$.
\changed{Let us observe}
 that $I_N$ belongs to the dual cone $\Pcone_\Vsp^*$ of $\Pcone_\Vsp$,
 that is,
$$
 0 < \cpling{y}{I_N} = y_{11} + \dots + y_{rr} \quad (y \in \lbar{\Pcone}_\Vsp \setminus \{0\}).
$$
\changed{
Indeed,
 since each $y_{kk}$ is a diagonal entry of the non-negative matrix
 $y \in \lbar{\Pcone}_\Vsp \setminus \{0\}$,
 we have $y_{kk} \ge 0$.
Suppose $\sum_{k=1}^r y_{kk} = 0$. 
Then $y_{11} = \dots = y_{rr} =0$ and $\mathrm{tr}\,y = \sum_{k=1}^r n_k y_{kk} =0$.
This together with the non-negativity of $y$ implies $y=0$, 
 which is a contradiction.
}  

For $T \in H_\Vsp$,
 define $\rho^*(T) \in GL(\Zsp_\Vsp)$ by
 $\cpling{y}{\rho^*(T)\eta} = \cpling{\rho(T)y}{\eta}\,\,\,(y,\,\eta \in \Zsp_\Vsp)$.
By \cite{V}, we have
 $\Pcone_\Vsp^* = \set{\rho^*(T)I_N}{T \in H_\Vsp}$.
Moreover,
 the map $H_\Vsp \owns T \mapsto \rho^*(T)I_N \in \Pcone_\Vsp^*$ is a diffeomorphism.
For $\usigma = (\sigma_1, \dots, \sigma_r) \in \complex^r$,
 we define the one-dimensional representation
 $\chi_{\usigma} : H_\Vsp \to \complex^{\times}$ by
 $\chi_{\usigma}(T) := (t_{11})^{2 \sigma_1} \dots (t_{rr})^{2 \sigma_r}\,\,\,(T \in H_\Vsp)$.
Note that any one-dimensional representation $\chi$ of $H_{\Vsp}$ is of the form $\chi_{\usigma}$,
 so that $\chi$ is determined by the values on the subgroup
 $A_\Vsp \subset H_\Vsp$ consisting of diagonal matrices.

Let us give some examples. 
When $n_1 = n_2 = \cdots = n_r =1$ and $\Vsp_{lk}=\real$ for all $1 \le k < l \le r$,
 the conditions (V1)--(V3) are clearly satisfied,
 and we have
 $\Zsp_\Vsp=\mathrm{Sym}(r,\real)$ and $\Pcone_\Vsp = \Pi_r$. 
For the case $r=3$, $n_1 = n_2 = n_3 = 1$, $\Vsp_{21} =\{0\}$ and $\Vsp_{31}=\Vsp_{32} =\real$,
 the space $\Zsp_\Vsp$ equals $\Zsp$ in (\ref{eqn:dual-Vinberg}).  
 
Let us set $r=3$, $n_1=2,\,n_2 = n_3 = 1$,
$$
\Vsp_{21} = \set{\begin{pmatrix} v & 0 \end{pmatrix}}{v \in \real},\quad
\Vsp_{31} = \set{\begin{pmatrix} 0 & v \end{pmatrix}}{v \in \real},
$$  
and $\Vsp_{32} = \{0\}$. 
Then we have
\begin{equation} \label{eqn:Vinberg_sp}
\Zsp_\Vsp = \set{
\begin{pmatrix} y_{11} & 0 & y_{21} & 0\\ 0 & y_{11} &0 & y_{31} \\ 
 y_{21} & 0 & y_{22} & 0 \\ 0 & y_{31} & 0 & y_{33} \end{pmatrix}}
 {y_{11},\,y_{22},\,y_{33},\,y_{21},\,y_{31} \in \real}
\end{equation}
and 

\begin{equation} \label{eqn:Vinberg_cone} 
\begin{aligned}
\Pcone_{\Vsp} &= \set{y \in \Zsp_\Vsp}{y \mbox{ is positive definite}} \\
&= \set{y \in \Zsp_\Vsp}{y_{11} >0,\,\,y_{11} y_{22} - (y_{21})^2 >0,\,\,y_{11}y_{33} - (y_{31})^2>0},
\end{aligned}
\end{equation}
 which is exactly the Vinberg cone \cite{V}.

Set $r=2,\,\,n_1 = m \ge 1,\, n_2 =1$ 
 and $\Vsp_{21} = \mathrm{Mat}(1,m;\real)$.
Then
\begin{align}
\Zsp_\Vsp &= \set{
\begin{pmatrix} y_{11}& & & v_1 \\ & \ddots & & \vdots\\ & & y_{11} & v_m \\ v_1& \cdots & v_m & y_{22}
\end{pmatrix}}{y_{11},\,y_{22},\,v_1,\dots, v_m \in \real}, \nonumber \\
\Pcone_\Vsp &= \set{y \in \Zsp_\Vsp}{y_{11} >0,\,\,y_{11} y_{22} - (v_1)^2 - \dots - (v_m)^2 >0},
\label{eqn:Lorentz}
\end{align} 
so that we obtain the Lorentz cone of dimension $m+2$.\\

\noindent{\bf 3.3. ~Basic quadratic maps.}\indent
Let $W_\Vsp^i\,\,\,(i=1, \dots, r)$ be 
 the subspace of $\mathrm{Mat}(N, n_i;\real)$ consisting of
 matrices $x$ of the form
$$
 x = \begin{pmatrix} 0_{n_1 + \dots + n_{i-1}, n_i} \\ X_{ii}
   \\ \vdots \\ X_{ri} \end{pmatrix}
 \quad 
 \left(
 \begin{aligned}
 X_{ii} &= x_{ii} I_{n_i},\,\,x_{ii} \in \real \\
 X_{li} &\in \Vsp_{li} \,\,\,(l=i+1, \dots, r)
 \end{aligned}
 \right).
$$
For example,
 when $\Zsp_{\Vsp}$ is the one in (\ref{eqn:Vinberg_sp}), 
 we have
\begin{align*} 
 W_\Vsp^1 
 &= \set{\begin{pmatrix} x_{11} & 0 \\ 0 & x_{11} \\ x_{21} & 0 \\ 0 & x_{31} \end{pmatrix}}
  {x_{11},\,x_{21},\,x_{31} \in \real},\\
 W_\Vsp^2 
 &= \set{\begin{pmatrix} 0 \\ 0 \\ x_{22} \\ 0 \end{pmatrix}}{x_{22} \in \real},\qquad
 W_\Vsp^3 
 = \set{\begin{pmatrix} 0 \\ 0 \\ 0 \\ x_{33} \end{pmatrix}}{x_{33} \in \real},
\end{align*}
 while
 for the case that $\Zsp_\Vsp$ is the space $\Zsp$ in (\ref{eqn:dual-Vinberg}),
 we have
\begin{align*}
 W_\Vsp^1 &= \set{\begin{pmatrix} x_{11} \\ 0 \\ x_{31}\end{pmatrix}}{x_{11},\,x_{31} \in \real},\quad
 W_\Vsp^2 = \set{\begin{pmatrix} 0 \\ x_{22} \\ x_{32}\end{pmatrix}}{x_{22},\,x_{32} \in \real},\\
 W_\Vsp^3 &= \set{\begin{pmatrix} 0 \\ 0 \\ x_{33}\end{pmatrix}}{x_{33} \in \real}.
\end{align*}
For $T \in H_\Vsp$ and $x \in W_\Vsp^i$,
 we see from (V1) that  $Tx \in W_\Vsp^i$,
 which defines a representation $\tau_i :H_\Vsp \to GL(W_\Vsp^i)$. 
Since $W_\Vsp^i \subset \mathrm{Mat}(N,n_i;\real)$,
 we can consider the restriction $q_\Vsp^i$ of the $\Pi_N$-positive quadratic map
 $q_{N,n_i} : \mathrm{Mat}(N,n_i;\real) \to \mathrm{Sym}(N,\real)$ in Example 2
 to the space $W_\Vsp^i$.
Thanks to (V2) and (V3), 
 we have $q_{\Vsp}^i(x) = x \transp{x} \in \Zsp_\Vsp$
 for $x \in W_\Vsp^i$.
Then the quadratic map
 $q_\Vsp^i : W_\Vsp^i \owns x \mapsto x \transp{x} \in \Zsp_\Vsp$
 is $\Pcone_\Vsp$-positive.
On the other hand,
 we observe
\begin{equation} \label{eqn:q-taui}
 q_\Vsp^i(\tau_i(T)x) = (Tx)\transp{(Tx)} =\rho(T) q_{\Vsp}^i(x)
 \qquad (x \in W_\Vsp^i,\,\,T \in H_\Vsp), 
\end{equation}
 which yields the group homomorphism 
 $H_\Vsp \owns T \mapsto (\rho(T),\tau_i(T)) \in \mathrm{Aut}(\Pcone_\Vsp,q_\Vsp^i)$.
It follows that the quadratic map $q_\Vsp^i$ is homogeneous.
We call $q_\Vsp^1, \dots, q_\Vsp^r$ \textit{basic quadratic maps} for $\Pcone_\Vsp$.
Recalling the inner product on $\Vsp_{lk}$ given by (\ref{eqn:innerprod_Vlk}),
 we define an inner product on the space $W_\Vsp^i$ 
 via the natural isomorphism
\begin{equation} \label{eqn:decomp_WVi}
 W_\Vsp^i \simeq \real \oplus \dirsum_{l>i} \Vsp_{li}.
\end{equation}
Taking an orthonormal basis of $W_\Vsp^i$ with respect to 
 the inner product,
 we identify $W_\Vsp^i$ with $\real^{m(i)}$,
 where $m(i) := \dim W_\Vsp^i$.
Then we consider 
 the linear map $\phi_{\Vsp}^i : \changed{\Zsp_\Vsp \equiv \Zsp_\Vsp^*} \to \mathrm{Sym}(m(i), \real)$
 associated to the quadratic map $q_\Vsp^i$.
Note that,
 if we write $n_{li}$ for $\dim \Vsp_{li}$, we have $m(i) = 1 + \sum_{l>i} n_{li}$.


\begin{prop} \label{prop:dualcone}
{\rm (i)}
One has
 $$
 \det \tau_i(T) = \chi_{\um(i)/2}(T) \qquad (T \in H_\Vsp), 
 $$
 where 
 $\um(i) = (0,\cdots,0,1, n_{i+1,i},\dots, n_{ri}) \in \integer^r$.\\
{\rm (ii)}
For $\eta = \rho(T)^*I_N \in \Pcone_{\Vsp}^*$ with $T \in H_{\Vsp}$,
 one has
 $$
 \det \phi_\Vsp^i(\eta) = \chi_{\um(i)}(T).
 $$
In particular, 
 $(t_{rr})^2 = \det \phi_\Vsp^r(\eta)$.\\
{\rm (iii)}
For $1 \le i < r$,
 there exist integers $c_{i+1,i},\dots, c_{ri}$ such that
\begin{align*}
 (t_{ii})^2 
 &= \det \phi_\Vsp^i(\eta) \cdot 
 (\det \phi_{\Vsp}^{i+1}(\eta))^{c_{i+1,i}}\cdots (\det \phi_\Vsp^r(\eta))^{c_{ri}}\\
 & \qquad (\eta = \rho^*(T)I_N \in \Pcone_{\Vsp}^*,\,\,T \in H_{\Vsp}).
\end{align*}
{\rm (iv)}
One has
 $\Pcone_\Vsp^* = \set{\eta \in \changed{\Zsp_\Vsp}}{\det \phi_\Vsp^i(\eta) >0\,\,\,(i=1, \dots, r)}$.
\end{prop}

\pf
(i) 
Since
 $H_\Vsp \owns T \mapsto \det \tau_i(T) \in \complex^{\times}$
 is a one-dimensional representation,
 it is sufficient to check the equality for diagonal matrices $T \in A_\Vsp$. 
In this case,
 the isomorphism (\ref{eqn:decomp_WVi}) 
 gives the eigenspace decomposition of $\tau_i(T)$,
 where $\real$ and $\Vsp_{li}$
 correspond to the eigenvalues $t_{ii}$ and $t_{ll}$ respectively.
Therefore we have
 $\det \tau_i(T) = t_{ii} \prod_{l>i} t_{li}^{n_{li}} = \chi_{\um(i)/2}(T)$.\\
(ii) 
Thanks to (\ref{eqn:innerprod_Vlk}), 
 we have
 $(x|x) = \cpling{q_\Vsp^i(x)}{I_N}$,
 which implies $\phi_\Vsp^i(I_N) = I_{m(i)}$. 
Thus we get
 $\det \phi_\Vsp^i(\eta) = \det \phi_\Vsp^i(\rho(T)^*I_N) = (\det \tau_i(T))^2 
 = \chi_{\um(i)}(T)$
 by (\ref{eqn:dualg}) and (i).\\
(iii) 
We see from (ii) that 
$(t_{ii})^2 = \det \phi_\Vsp^i(\eta) \cdot (t_{i+1,i+1})^{-2 \changed{n_{i+1, i}}} \cdots (t_{rr})^{-2 \changed{n_{r i}}}$,
 whence we can deduce (iii) recursively.\\
(iv) 
It is known (\cite[Chapter 3, Section 3]{V1} and \cite[Section 1]{G1})
 that a homogeneous cone is described as the subset of the ambient vector space
 consisting of points 
 \changed{at which} all relative\changed{ly} invariant \changed{(appropriately normalized)} functions are positive. 
Thus the assertion follows from (iii).
\qed
$ $\\

\noindent
{\bf Example 4.}
Let $\Pcone_\Vsp$ be the Vinberg cone, that is, 
 $\Zsp_\Vsp$ is as in (\ref{eqn:Vinberg_sp}).
Then we have
$$
\phi_\Vsp^1(\eta) 
= \begin{pmatrix} \eta_{11} & \eta_{21} & \eta_{31} \\ \eta_{21} & \eta_{22} & 0 \\ 
  \eta_{31} & 0 & \eta_{33} \end{pmatrix},\quad
\phi_\Vsp^2(\eta)
= \eta_{22}, \quad
 \phi_{\Vsp}^3(\eta) = \eta_{33}
$$
 for $\eta \in \changed{\Zsp_\Vsp}$.
If $\eta = \rho^*(T)I_4$ for $T \in H_\Vsp$,
 we have
$$
 (t_{11})^2 (t_{22})^2(t_{33})^2 = \det \phi_\Vsp^1(\eta),\quad 
 (t_{22})^2 = \eta_{22}, \quad
 (t_{33})^2 = \eta_{33},
$$ 
 so that
$$
 (t_{11})^2 = \frac{\det \phi_\Vsp^1(\eta)}{\eta_{22} \eta_{33}} 
 = \eta_{11} - \frac{(\eta_{21})^2}{\eta_{22}} - \frac{(\eta_{31})^2}{\eta_{33}}.
$$
On the other hand, 
 we have by Proposition~\ref{prop:dualcone} (iv)
$$
\Pcone_\Vsp^* = \set{\eta \in \Zsp_\Vsp}
{\eta_{11} \eta_{22} \eta_{33} - \eta_{33}(\eta_{21})^2 - \eta_{22}(\eta_{31})^2 >0,\,\,
 \eta_{22} >0,\,\,\eta_{33} >0}.
$$
Therefore,
 if \changed{$\Zsp$ is the space in (\ref{eqn:dual-Vinberg}),}
 the linear isomorphism
$$
\changed{\iota:}
\Zsp_\Vsp \owns 
\begin{pmatrix} \eta_{11} & 0 & \eta_{21} & 0\\ 0 & \eta_{11} &0 & \eta_{31} \\ 
 \eta_{21} & 0 & \eta_{22} & 0 \\ 0 & \eta_{31} & 0 & \eta_{33} \end{pmatrix}
\mapsto
 \begin{pmatrix} \eta_{33} & 0 & \eta_{31} \\ 0 & \eta_{22} & \eta_{21}\\
                         \eta_{31} & \eta_{21} & \eta_{11} \end{pmatrix} 
\in \changed{\Zsp}
$$
 gives a bijection from $\Pcone_\Vsp^*$ onto $\Pcone = \Zsp \cap \Pi_3$.
\changed{
The adjoint map $\iota^* : \Zsp^* \to \Zsp_\Vsp^* \equiv \Zsp_\Vsp$
 gives the matrix realization of the Vinberg cone $\Qcone$
 as the homogeneous cone $\Pcone_\Vsp$. }
\\

\noindent{\bf 3.4. ~Standard quadratic maps and $H$-orbits in $\lbar{\Pcone_\Vsp}$.}\indent
We call the maps $q_\Vsp^{\uep}\,\,\,(\uep \in \{0,1\}^r,\,\uep \ne (0, \dots, 0))$
 \textit{standard quadratic maps}.
They are of particular importance
 among the virtual quadratic maps
 $q_\Vsp^{\us} = (q_\Vsp^1)^{\oplus s_1} \oplus \cdots \oplus  (q_\Vsp^r)^{\oplus s_r}\,\,\,
 (\us = (s_1,\dots, s_r) \in \real^r)$.  
Let us denote by $I(\uep)$ the set $\set{1 \le i \le r}{\uep_i =1}$. 
We identify $q_\Vsp^{\uep}$ with
 a direct sum
 $\sum_{i \in I(\uep)}^\oplus q_\Vsp^i$
 on the space 
 $W_\Vsp^{\uep} := \sum_{i \in I(\uep)}^\oplus W_\Vsp^i$.
Recalling (\ref{eqn:decomp_WVi}),
 we have the isomorphism
 $W_\Vsp^{\uep} \simeq \sum_{i \in I(\uep)}^\oplus (\real \oplus \sum_{l>i}^\oplus \Vsp_{li})$,
 which enables us to describe
 $y = q_\Vsp^{\uep}(x)\in \Zsp_{\Vsp}\,\,\,(x \in W_\Vsp^{\uep})$
 as the matrix composed of the blocks
\begin{equation} \label{eqn:x-ylk}
Y_{lk} = \sum_{i \le k,\,\,i \in I(\uep)} X_{li} \transp{X_{ki}}
 \quad (1 \le k \le l \le r),
\end{equation}
 where $X_{ii} := x_{ii} I_{n_i}$ for $i \in I(\uep)$. 
For the case $l=k$, 
 we have $Y_{kk} = y_{kk} I_{n_k}$ and 
\begin{equation} \label{eqn:x-ykk}
y_{kk} = \sum_{i \le k,\,\,i \in I(\uep)} \norm{X_{ki}}^2 
\end{equation}
 thanks to (\ref{eqn:innerprod_Vlk}), 
 where we put
 $\norm{X_{ii}} := |x_{ii}|$.
\changed{
For each $x \in W_\Vsp^{\uep}$,
 let $T_x \in \mathrm{Mat}(N, \real)$ be a lower triangular matrix 
 whose $(k,i)$-block component is $X_{ki}$ for $k>i$ with $i \in I(\ep)$,
 and other components are zero.
Then we have 
\begin{equation} \label{eqn:Tx}
 q_{\Vsp}^{\uep}(x) = T_x \,\transp{T_x}.
\end{equation}
For example,
 if $r=3$ and $\uep = (1,0,1)$,
 then an element $x$ of $W_{\Vsp}^{\uep} = W_{\Vsp}^1 \oplus W_{\Vsp}^3$ is of the form
$$
 x = \begin{pmatrix} X_{11} \\ X_{21} \\ X_{31}\end{pmatrix}
 \oplus
 \begin{pmatrix} 0 \\ 0 \\ X_{33} \end{pmatrix}
 \qquad
 \left(
 \begin{aligned} {} & X_{11} = x_{11} I_{n_1},\,\,X_{33} = x_{33} I_{n_3}\\
  {} & x_{11},\,x_{33}\in \real,\,\,X_{21} \in \Vsp_{21},\,X_{31} \in \Vsp_{31}
 \end{aligned}
 \right),
$$   
 and we have
$$
 T_x = \begin{pmatrix} X_{11} &  & \\ X_{21} &0 & \\ X_{31} & 0 & X_{33} \end{pmatrix} \in \mathrm{Mat}(N, \real).
$$
}

For $\uep \in \{0,1\}^r$,
 let $E_{\uep}$ be 
 the element of $\Vsp$ given by
$$
 E_{\uep} := \begin{pmatrix} \ep_1 I_{n_1} & & \\ & \ddots & \\ & & \ep_r I_{n_r} \end{pmatrix},
$$ 
 and $\orbt_{\uep}$ the $H_\Vsp$-orbit
 $\rho(H_{\Vsp})E_{\uep} \subset \Zsp_\Vsp$
 through $E_{\uep}$.
In particular, 
 the orbit
 $\orbt_{(0,\dots,0)}$ is the origin $\{0\}$, 
 while $\orbt_{(1,\dots,1)} = \rho(H_\Vsp)I_N$ equals the cone $\Pcone_{\Vsp}$. 
It is shown in \cite[Theorem 3.5]{I1} that the $H_\Vsp$-orbit decomposition
 of the closure $\lbar{\Pcone_\Vsp}$ is given as 
$$
\lbar{\Pcone_\Vsp} = \bigsqcup_{\uep \in \{0,1\}^r} \orbt_{\uep}.
$$


\begin{prop} \label{prop:orbit}
If $\uep \ne (0,\dots, 0)$, 
 the image of the quadratic map $q_\Vsp^{\uep}$ equals 
 the closure $\lbar{\orbt_{\uep}}$ of the orbit $\orbt_{\uep}$.
\end{prop}

\pf
\changed{
Let $W_{\Vsp}^{\uep, +}$ be the subset 
 $\set{x \in W_\Vsp^{\uep}}{x_{ii}>0\,\,(i \in I(\uep))}$
 of $W_\Vsp^{\uep}$.
For $y = \rho(T)E_{\uep} = T E_{\uep}\,\transp{T} \in \orbt_{\uep}$
 with $T \in H_\Vsp$,
 we take a unique $x \in W_{\Vsp}^{\uep, +}$ for which $T_x$ equals $TE_{\uep}$.
Then we have $y = q_\Vsp^{\uep}(x)$ by (\ref{eqn:Tx}).
Conversely, for any $x \in W_{\Vsp}^{\uep,+}$,
 we put $T := T_x + I - E_{\uep} \in H_{\Vsp}$
 so that 
 $q_\Vsp^{\uep}(x) = T E_{\uep} \,\transp{T} \in \orbt_{\uep}$.
Therefore we obtain
 $\orbt_{\uep} = q_\Vsp^{\uep}(W_{\Vsp}^{\uep,+})$.
}
On the other hand,
 putting
 $W_{\Vsp}^{\uep, *} := \set{x \in W_\Vsp^{\uep}}{x_{ii} \ne 0 \,\,\,(i \in I(\uep))}$,
 we see easily that
 $q_\Vsp^{\uep}(W_{\Vsp}^{\uep,*}) = q_\Vsp^{\uep}(W_{\Vsp}^{\uep,+})$.
Since $W_{\Vsp}^{\uep, *}$ is an open dense subset of $W_\Vsp^{\uep}$, 
 the orbit $\orbt_{\uep}$ is dense in the image of the quadratic map $q_\Vsp^{\uep}$,
 which is necessarily closed.
Indeed, 
 introducing the projective imbedding $\iota_V$ of a vector space $V$ by 
 $\iota_V: V \owns y \mapsto [1,y] \in \mathbf{P}V := (\real \times V \setminus \{(0,0)\}) / \changed{\real^{\times}}$,
 we can extend
 $q_\Vsp^{\uep}: W_\Vsp^{\uep} \to \Zsp_\Vsp$
 to the map  
 $\tilde{q}_\Vsp^{\uep} : \mathbf{P}W_\Vsp^{\uep} \owns [t,x] \mapsto [t^2,q_\Vsp^{\uep}(x)] \in \mathbf{P}\Zsp_\Vsp$
 because
 $\changed{q_\Vsp^{\uep}}(x) \ne 0$ for $x\ne 0$.
The image $\tilde{q}_\Vsp^{\uep}(\mathbf{P}W_\Vsp^{\uep})$ is compact,
 so that
 $q_\Vsp^{\uep}(W_\Vsp^{\uep}) = \iota_{\Zsp_\Vsp}^{-1}(\tilde{q}_\Vsp^{\uep}(\mathbf{P}W_\Vsp^{\uep}))$ is closed.
%
%
\qed

\begin{rem}
In the proof of Proposition \ref{prop:orbit}, 
 we see that the quadratic map
 $q_\Vsp^{\uep}$ gives a surjective map from $W_\Vsp^{\uep, +}$ onto $\orbt_{\uep}$.
The map is also one-to-one thanks to \cite[Lemma 3.3 (ii)]{I1},
 while 
 the map $q_\Vsp^{\uep}: W_\Vsp^{\uep, *} \to \orbt_{\uep}$ is $2^{\sharp I(\uep)}$-to-one.
Actually, 
 a large part of the content of this section is presented in language of normal $j$-algebra
 in \cite[Sections 3 and 4]{I1}.
 \end{rem}

We define the representation $\tau_{\uep} : H_\Vsp \to GL(W_\Vsp^{\uep})$
 as the direct sum of the representations
 $(\tau_i,\, W_\Vsp^i)$ of $H_{\Vsp}$ for $i \in I(\uep)$. 
Then we have by (\ref{eqn:q-taui})
\begin{equation} \label{eqn:q-tauep}
 q_\Vsp^{\uep}(\tau_{\uep}(T)x) = \rho(T) q_\Vsp^{\uep}(x) 
 \qquad (x \in W_\Vsp^{\uep},\,\,T \in H_\Vsp),
\end{equation}
 which implies that $q_\Vsp^{\uep}$ is homogeneous.

The open set $W_\Vsp^{\uep, +} \subset W_\Vsp^{\uep}$ is preserved
 by the action of $\tau_{\uep}(H_\Vsp)$.
We put
$$
R_+(\uep) 
 := \set{\uu = (u_1, \dots, u_r) \in \real^r}
        {u_i = 0\,\,(\mbox{if }\ep_i = 0),\quad u_i>0\,\,(\mbox{if }\ep_i=1)}.
$$
For $\uu \in R_+(\uep)$,
 let $\meas_{\uu}^{\uep}$ be the measure on $W_\Vsp^{\uep, +}$
 given by
\begin{align} 
 \meas_{\uu}^{\uep}(dx) 
 & := \prod_{i \in I(\uep)} \Bigl\{ 
 \frac{ 2 (x_{ii})^{2 u_i - 1} \,dx_{ii} }{\Gamma(u_i)} \cdot 
 \prod_{l>i} \frac{dX_{li}}{\pi^{n_{li}/2}}
 \Bigr\} \nonumber\\
 &= \frac{\prod_{i \in I(\uep)} (x_{ii})^{2 u_i - 1} }{\Gamma_{\uep}(\uu)}\, dx
 \label{eqn:def_of_meas_u}
 \qquad (x \in W_\Vsp^{\ep, +}),
\end{align}
 where 
 $\Gamma_{\uep}(\uu) = \pi^{\dim W_\Vsp^{\uep} /2} \prod_{i \in I(\uep)} \frac{\Gamma(u_i)}{2 \sqrt{\pi}}$.
When $\uu = \uep/2$,
 the measure $\meas_{\uep/2}^{\uep}$ equals a constant multiple of the Lebesgue measure,
 that is,
 $\meas_{\uep/2}^{\uep}(dx) = 2^{\sharp I(\uep)} \pi^{-\dim W_\Vsp^{\uep} /2} dx$. 
We define $\up(\uep) := (p_1(\uep),p_2(\uep),\dots, p_r(\uep))$ by
 $$ p_k(\uep) := \sum_{l>k} \ep_l n_{lk}. $$


\begin{lemma} \label{lemma:meas_u}
{\rm (i)}
For a measurable set $A \subset W_\Vsp^{\uep, +}$,
 one has 
$$
 \meas_{\uu}^{\uep}(\tau_{\uep}(T)A) = \chi_{\uu + \up(\uep)/2}(T) \meas_{\uu}^{\uep}(A)
 \qquad (T \in H_{\Vsp}).
$$
{\rm (ii)}
One has
\begin{equation} \label{eqn:equal1}
\int_{W_\Vsp^{\uep,+}} e^{-\norm{x}^2}\,\meas_{\uu}^{\uep}(dx) = 1.
\end{equation}
\end{lemma}

\pf (i)
If $x' = \tau_{\uep}(T)x \in W_\Vsp^{\uep, +}$ with $x \in W_\Vsp^{\uep,+}$,
 we have $x'_{ii} = t_{ii} x_{ii}$ for $i \in I(\uep)$.
Thus
\begin{equation} \label{eqn:chi_u}
 \prod_{i\in I(\uep)} (x'_{ii})^{2 u_i-1} 
= \chi_{\uu - \uep/2}(T) \prod_{i\in I(\uep)} (x_{ii})^{2 u_i-1}. 
\end{equation}
On the other hand,
 we observe that
$$
 dx' = |\det \tau_{\uep}(T)| dx = \Bigl(\prod_{i \in I(\uep)} \det \tau_i(T) \Bigr) dx,
$$
 and the last term equals 
$\Bigl(\prod_{i \in I(\uep)} \chi_{\um(i)/2}(T) \Bigr) dx$
 by Proposition~\ref{prop:dualcone} (i).
Since
$$
 \sum_{i \in I(\ep)} \um(i) = \sum_{i=1}^r \ep_i \um(i) = \uep + \up(\uep),
$$
 we have
$dx' = \chi_{\uep/2 + \up(\uep)/2}(T) dx$,
which together with (\ref{eqn:chi_u}) implies (i).\\
(ii)
By definition, we have
$$
\norm{x}^2 = \sum_{i \in I(\uep)} \Bigl\{(x_{ii})^2 + \sum_{l>i} \norm{X_{li}}^2 \Bigr\}.
$$
Thus the left-hand side of (\ref{eqn:equal1}) equals
$$
\prod_{i \in I(\uep)} 
\Bigl\{
\frac{2}{\Gamma(u_i)}\int_0^{\infty} e^{-(x_{ii})^2} (x_{ii})^{2u_i-1} dx_{ii} 
\prod_{l>i} \int_{\Vsp_{li}} e^{-\norm{X_{li}}^2} \frac{dX_{li}}{\pi^{n_{li}/2}}
\Bigr\}.
$$
Therefore
 we obtain (\ref{eqn:equal1}) from 
 $\int_{\Vsp_{li}} e^{-\norm{X_{li}}^2} dX_{li} = \pi^{n_{li}/2}$
 and
 $\int_0^{\infty} e^{-(x_{ii})^2} (x_{ii})^{2u_i-1} dx_{ii} = \Gamma(u_i)/2$.
\qed
$ $\\
\begin{rem}\label{rem-aboutM}  
 Lemma~\ref{lemma:meas_u} tells us that 
 $e^{-\norm{x}^2}\,\meas_{\uu}^{\uep}(dx)$
 is a probability measure on $W_{\Vsp}^{\uep,+}$.
Actually, we see from the proof that 
 if $X^{\uu}$ is an $W_{\Vsp}^{\uep,+}$-valued random variable
 with the law $e^{-\norm{x}^2}\,\meas_{\uu}^{\uep}(dx)$,
 then its components \changed{are independent and }
 satisfy $\sqrt{2} X^{\uu}_{li} \sim N(0,I_{n_{li}})$,
 and $(\sqrt{2}X^{\uu}_{ii})^2 \sim \chi^2(2u_i)$,
 where
 $\chi^2(u)$ denotes the chi-square law with the density
 $2^{-u}\Gamma(u/2)^{-1}e^{-t/2} t^{u-1}\,\,\,(t>0)$.
We shall see later that any Wishart law associated to a homogeneous quadratic map
 is the image of this measure by an appropriate quadratic map.\\
\end{rem}
\noindent{\bf 3.5. ~Gindikin-Riesz distributions.}\indent
For $\usigma = (\sigma_1, \dots, \sigma_r) \in \complex^r$,
 we denote $(\sigma_r, \dots, \sigma_1)$ by $\usigma^*$.
Let $\Delta^*_{\usigma}$ be the function on the cone $\Pcone_{\Vsp}^*$
 given by
\begin{equation} \label{eqn:def_of_Deltastar}
 \Delta^*_{\usigma}(\rho^*(T)I_N) = \chi_{\usigma^*}(T) \qquad (T \in H_\Vsp). 
\end{equation}
By Proposition~\ref{prop:dualcone} (ii) and (iii),
 $\Delta^*_{\usigma}(\eta)$ can be expressed as a product of powers of 
 the polynomials $\det \phi_\Vsp^i(\eta)$.
Putting
$$
 E_{\underline{\eta}} := \begin{pmatrix} \eta_1 I_{n_1} & & \\ & \ddots & \\ & & \eta_r I_{n_r} \end{pmatrix} \in \Zsp_\Vsp
$$
 for $\underline{\eta} = (\eta_1, \dots, \eta_r) \in \real^r_{>0}$,
 we have
\begin{equation}\label{eqn:DeltastarEeta}
 \Delta^*_{\usigma}(E_{\underline{\eta}}) = (\eta_1)^{\sigma_r} (\eta_2)^{\sigma_{r-1}} \dots (\eta_r)^{\sigma_1}.
\end{equation}
For each $\usigma \in \complex^r$,
Gindikin (\cite{G1}, \cite{G2}) constructed a tempered distribution
 $\Rz_{\usigma} \in \mathcal{S}'(\Zsp_\Vsp)$ 
 whose Laplace transform $L_{\Rz_{\usigma}}(\theta) = \Rz_{\usigma}(e^{\cpling{y}{\theta}})$ is given by
\begin{equation} \label{eqn:GR-Laplace}
 L_{\Rz_{\usigma}}(\theta) = \Delta^*_{-\usigma^*}(-\theta) \qquad (\theta \in  -\Pcone_\Vsp^*).
\end{equation}
We call $\Rz_{\usigma}$
 the \textit{Gindikin-Riesz distribution} on the homogeneous cone $\Pcone_\Vsp$.
The support of $\Rz_{\usigma}$ is contained in $\lbar{\Pcone_\Vsp}$,
 and $\Rz_{\usigma}$ is relatively invariant under the action of $\rho(H_\Vsp)$,
 that is,
\changed{
\begin{equation} \label{eqn:GR_T}
 \Rz_{\usigma}(f \circ \rho(T)) = \chi_{-\usigma}(T)\Rz_{\usigma}(f)
\end{equation}
} for $T \in H_\Vsp$ and $f \in \mathcal{S}(\Zsp_\Vsp)$.


\begin{prop} \label{prop:image_meas}
For non-zero $\uep \in \{0,1\}^r$ and $\uu \in R_+(\uep)$, 
 put $\usigma := \uu + \up(\uep)/2$.
Then $\Rz_{\usigma}$
 is the image of $\meas_{\uu}^{\uep}$ by the standard quadratic map $q_\Vsp^{\uep}$.   
\end{prop}

\pf
By (\ref{eqn:GR-Laplace}),
 it is sufficient to show that
$$
 \int_{W_\Vsp^{\uep,+}} e^{\cpling{q_{\Vsp}^{\uep}(x)}{\theta}}\,\meas_{\uu}^{\uep}(dx) 
 = \Delta^*_{-\usigma^*}(-\theta)
$$
 for $\theta \in -\Pcone_\Vsp^*$.
Take $T \in H_\Vsp$ for which $\theta = -\rho^*(T)I_N$.
Then the left-hand side is
$$
 \int_{W_\Vsp^{\uep,+}} e^{-\cpling{q_{\Vsp}^{\uep}(x)}{\rho^*(T)I_N}}\,\meas_{\uu}^{\uep}(dx)
 = \int_{W_\Vsp^{\uep,+}} e^{-\cpling{q_{\Vsp}^{\uep}(\tau_{\uep}(T)x)}{I_N}}\,\meas_{\uu}^{\uep}(dx)
$$
 by (\ref{eqn:q-tauep}), 
 and it is equal to
 $$ \chi_{-(\uu + \up(\uep)/2)}(T)
 \int_{W_\Vsp^{\uep,+}} e^{-\cpling{q_{\Vsp}^{\uep}(x)}{I_N}}\,\meas_{\uu}^{\uep}(dx) $$
 by Lemma \ref{lemma:meas_u} (i).
Since $\cpling{q_{\Vsp}^{\uep}(x)}{I_N} = \norm{x}^2$ by (\ref{eqn:x-ykk}),
 we see from Lemma \ref{lemma:meas_u} (ii) that
$$
 \int_{W_\Vsp^{\uep,+}} e^{\cpling{q_{\Vsp}^{\uep}(x)}{\theta}}\,\meas_{\uu}^{\uep}(dx)
 = \chi_{-(\uu + \up(\uep)/2)}(T) = \Delta^*_{-\usigma^*}(-\theta).
$$
\qed
$ $\\
 
We set
\begin{align*}
\Xi(\uep) &:= \set{\usigma = \uu + \up(\uep)/2}{\uu \in R_+(\uep)}\\
 &= \set{\usigma \in \real^r}
   {\sigma_i = p_i(\uep)/2 \,\,(\mbox{if }\ep_i = 0), \quad 
    \sigma_i > p_i(\uep)/2 \,\,(\mbox{if }\ep_i = 1)}.
\end{align*}
If $\uep \ne (0,\dots,0)$ and $\usigma \in \Xi(\uep)$, 
 then $\Rz_{\usigma}$ is a positive measure on the orbit $\orbt_{\uep}$
 by Proposition~\ref{prop:image_meas}.
For the case $\ep = (0,\dots, 0)$,
 we have $\Xi(0,\dots,0) = \{(0,\dots,0)\}$
 and $\Rz_{(0,\dots,0)}$ is the Dirac measure at the origin $\{0\}$.
It is proven in \cite{I1} that 
 they exhaust all the cases that 
 $\Rz_{\usigma}$ is a positive measure.


\begin{thm}[\mbox{\cite[Theorem 6.2]{I1}}] \label{thm:Gindikinset}
The Gindikin-Riesz distribution $\Rz_{\usigma}$ is a positive measure
 if and only if $\usigma \in \Xi := \bigsqcup_{\uep \in \{0,1\}^r} \Xi(\ep)$.
Moreover, if $\usigma \in \Xi(\uep)$, 
 then $\Rz_{\usigma}$ is a measure on $\orbt_{\uep}$.
\end{thm}


The parameter set $\Xi$ is also described as
 $$
 \Xi = \set{\sum_{i=1}^r \ep_i (0,\cdots,0,u_i, n_{i+1,i}/2,\dots, n_{ri}/2) }
 {\ep_i \in \{0,1\},\,u_i >0\,\,(i=1, \dots, r)}.
 $$

\noindent{\bf \changed{3.6}.~Riesz measures and Gindikin-Riesz distributions.} \indent
Let us investigate a relation of the Riesz measure\changed{s}
 $\mu_q$ associated to 
 \changed{homogeneous $\Pcone_{\Vsp}$-positive} 
 quadratic maps $q$ and
 the Gindikin-Riesz distributions on $\Pcone_\Vsp$. 


\begin{prop} \label{prop:muq-Riesz}
For $i=1, \dots, r$,
 one has $\mu_{q_\Vsp^i} = \pi^{m(i)/2} \Rz_{\um(i)}$.
\end{prop}
\pf
>From Lemma~\ref{lemma:Laplace_trans} and Proposition~\ref{prop:dualcone} (ii),
 we have
\begin{equation} \label{eqn:Laplace_muqVi}
 L_{\mu_{q_\Vsp^i}}(\theta) = \pi^{m(i)/2} \det \phi_\Vsp^i(-\theta)^{-1/2}
= \pi^{m(i)/2} \Delta^*_{-\um(i)^*/2}(-\theta)
 \qquad (\theta \in -\Pcone_\Vsp^*),
\end{equation}
 which implies the statement.
\qed
$ $\\
 
Assume that there exists the Riesz measure
 $\mu_{q_\Vsp^{\us}}$ 
 associated to a virtual quadratic map
 $q_\Vsp^{\us} = (q_\Vsp^1)^{\oplus s_1} \oplus \cdots \oplus  (q_\Vsp^r)^{\oplus s_r}$.
By (\ref{eqn:virtual_LaplaceTrans}) and (\ref{eqn:Laplace_muqVi}),
 we have for $\changed{\theta} \in -\Pcone_\Vsp^*$
$$
 L_{\mu_{q_\Vsp^{\us}}}(\changed{\theta}) = \prod_{i=1}^r \Bigl( \pi^{m(i)/2} \Delta^*_{-\um(i)^*/2}(-\changed{\theta}) \Bigr)^{s_i}.
$$
We put
\begin{equation} \label{eqn:brigde}
 \usigma := \sum_{i=1}^r s_i \um(i)/2
 = \frac{1}{2}\sum_{i=1}^r s_i (0, \dots, 0, 1, n_{i+1,i}, \cdots, n_{ri}).
\end{equation}
Then the equality above can be rewritten as
$$
 L_{\mu_{q_\Vsp^{\us}}}(\changed{\theta}) = \pi^{|\usigma|} \Delta^*_{- \usigma^*}(-\changed{\theta}),
$$
 where $|\usigma| := \sigma_1 + \dots + \sigma_r$. 
Thus $\mu_{q_\Vsp^{\us}}$ equals $\pi^{|\usigma|} \Rz_{\usigma}$, 
 so that
 $\usigma$ belongs to $\Xi(\uep)$ for some $\uep \in \{0,1\}^r$ 
 owing to Theorem~\ref{thm:Gindikinset}.
The converse argument is also valid. 
Therefore we obtain


\begin{thm} \label{thm:Riesz-Gindikin}
For a virtual quadratic map 
 $q_\Vsp^{\us} = (q_\Vsp^1)^{\oplus s_1} \oplus \cdots \oplus  (q_\Vsp^r)^{\oplus s_r}$,
 there exists the associated Riesz measure $\mu_{q_\Vsp^{\us}}$
 if and only if $\usigma := \sum_{i=1}^r s_i \um(i)/2$ belongs to $\Xi$.
In this case 
 $\mu_{q_\Vsp^{\us}} = \pi^{|\usigma|} \Rz_{\usigma}$,
 and there exist
 $\uep \in \{0,1\}^r$ and $\uu \in R_+(\uep)$ for which
 $\usigma = \uu + \up(\uep)/2$.
If $\uep \ne (0,\dots,0)$,
 $\mu_{q_\Vsp^{\us}}$ is the image of the measure 
 $\pi^{|\usigma|} \meas_{\uu}^{\uep}$ on $W_\Vsp^{\uep}$
 by the standard quadratic map $q_\Vsp^{\uep}$.
\end{thm}


Let $q: \real^m \to \Zsp_{\Vsp}$ be
 any \changed{homogeneous $\Pcone_\Vsp$-positive} quadratic map.
As is noted in Section 3.1,
 the group $pr_1(\mathrm{Aut}(\Pcone_{\Vsp},q))$
 acts on the cone $\Pcone_\Vsp$ transitively.
Assume first that $pr_1(\mathrm{Aut}(\Pcone_{\Vsp},q))$ contains $\rho(H_\Vsp)$.
Then the polynomial $\det \phi_q(\eta)$ is relatively invariant under the action of $\rho(H_\Vsp)$
 by (\ref{eqn:RelInv}).
Namely,
 for each $T \in H_\Vsp$,
 there exists $c_T >0$ such that
 $\det \phi_q(\rho^*(T)\eta) = c_T \det \phi_q(\eta)\,\,\,(\eta \in (\real^n)^*)$.  
It is easy to see that 
 the correspondence $H_\Vsp \owns T \mapsto c_T \in \real_{>0}$ is
 a one-dimensional representation, 
 so that
 we have $c_T = \chi_{\um}(T)$ for some $\um \in \real^r$.
Thus
 $\det \phi_q(\rho^*(T)I_n) =  \chi_{\um}(T) \det \phi_q(I_N)$
 for $T \in H_{\Vsp}$,
 which means that
 $\det \phi_q(\eta) = C \Delta^*_{\um^*}(\eta)$
 for $\eta \in \Pcone^*_\Vsp$,
 where $C := \det \phi_q(I_N)$.
By (\ref{eqn:DeltastarEeta}),
 we have
 $$ 
 \det \phi_q(E_{\underline{\eta}}) = C\, (\eta_1)^{m_1} \dots (\eta_r)^{m_r},
 $$ 
 which gives a practical way to determine $m_i$.
Indeed,
 we see from this formula that $m_i$ are non-negative integers. 
Comparing the degree\changed{s} of both sides, 
 we obtain $m= m_1 + \cdots + m_r$.
Similarly to Proposition~\ref{prop:muq-Riesz}, 
 we have
\begin{equation}
  \mu_q = C^{-1/2} \pi^{m/2} \Rz_{\um/2}.
\end{equation}
Let us \changed{consider} the virtual quadratic map $q^{\oplus s}$.
The associated Riesz measure exists if and only if $s\um/2 \in \Xi$,
 and in this case
$$
 \mu_{q^{\oplus s}} = C^{-s/2} \pi^{sm/2} \Rz_{s \um/2}.
$$
As for the general case,
 we have the following result.


\begin{prop} \label{prop:general_q}
Let $q: \real^m \to \Zsp_{\Vsp}$ be
 a \changed{homogeneous $\Pcone_\Vsp$-positive} quadratic map.
Then there exist $g_0 \in G(\Pcone_\Vsp)$, $\um \in \integer^r$,
 and $C >0$ for which
 $$
 \det \phi_q((g_0^{-1})^*\eta) = C \Delta^*_{\um^{\changed{*}}}(\eta) \quad (\eta \in \Pcone_\Vsp).
 $$
The Riesz measure $\mu_{q^{\changed{\oplus s}}}$ associated to the virtual quadratic map $q^{\changed{\oplus s}}$ exists
 if and only if $s \um /2 \in \Xi$.
In this case,
 $\mu_{q^{\changed{\oplus s}}}$ equals the image of $C^{-s/2} \pi^{s m/2}\Rz_{s \um/2}$ by $g_0$. 
\end{prop}

\pf
We note that $pr_1(\mathrm{Aut}(\Pcone_{\Vsp},q))$ acts on the cone $\changed{\Pcone_\Vsp}$ transitively,
 and that the identity component of $pr_1(\mathrm{Aut}(\Pcone_{\Vsp},q))$ equals the identity component of 
 an algebraic group (cf. \changed{\cite[Theorem 2]{I7}}).
It follows that
 an Iwasawa subgroup (maximal connected split solvable subgroup)
 $\mathcal{H}$ 
 \changed{of $pr_1(\mathrm{Aut}(\Pcone_{\Vsp},q))$}
 acts on $\Pcone_\Vsp$ simply transitively
 \changed{(\cite[Chapter 1]{V1})}.
Since $\mathcal{H}$ is also an Iwasawa subgroup of $G(\Pcone_\Vsp)$,
 it is conjugate to another Iwasawa subgroup $\rho(H_{\Vsp}) \subset G(\Pcone_\Vsp)$.
Namely,
 there exists $g_0 \in G(\Pcone_\Vsp)$ for which
 $g_0^{-1} \mathcal{H}g_0 = \rho(H_{\Vsp})$.
Let $q'$ be the $\Pcone$-positive quadratic map
  $g_0^{-1} \circ q : \real^m \to \Zsp_\Vsp$.
We have $\phi_{q'}(\eta) = \phi_q((g_0^{-1})^*\eta)$ for $\eta \in (\real^n)^*$
 because
 $\transp{x}\phi_{q'}(\eta)x = \cpling{q'(x)}{\eta} = \cpling{q(x)}{(g_0^{-1})^*\eta} 
 = \transp{x}\phi_{q}((g_0^{-1})^*\eta)x$
 for $x \in \real^m$.
It is easy to see that
$$
 \mathrm{Aut}(\Pcone_{\Vsp}, q') 
 =\set{(g_0^{-1} g_1 g_0, g_2) \in GL(\Zsp_\Vsp) \times GL(\real^m)}{(g_1,g_2) \in \mathrm{Aut}(\Pcone_{\Vsp},q)}.
$$
Then 
 $pr_1(\mathrm{Aut}(\Pcone_{\Vsp}, q'))
 =g_0^{-1}  pr_1(\mathrm{Aut}(\Pcone_{\Vsp}, q)) g_0
 \changed{\supset} g_0^{-1} \mathcal{H}g_0 = \rho(H_{\Vsp})$.
Thus we can apply the argument preceding Proposition~\ref{prop:general_q}
 \changed{for $q'$},
 so that we have
 $$
 \det \phi_q((g_0^{-1})^*\eta) = \det \phi_{q'}(\eta) = C \Delta^{\changed{*}}_{\um^{\changed{*}}}(\eta) 
 $$
 with some $C>0$ and $\um \in \integer^r$.
\changed{
Moreover $\mu_{(q')^{\oplus s}}$ equals $C^{-s/2} \pi^{s m/2}\Rz_{s \um/2}$
 if $s \um /2 \in \Xi$.
Since $q^{\oplus \pg{s}} = g_{\pg{0}} \circ (q')^{\oplus s}$,
 we  \pg{get} the last statement from Proposition~\ref{prop:virtual_trans}.}
\qed
$ $\\

Proposition~\ref{prop:general_q} states that the Riesz measure $\mu_q$ associated to a homogeneous $q$
 is equal to some Gindikin-Riesz distribution up to a linear transform on $G(\Pcone_\Vsp)$.
On the other hand,
 Theorem~\ref{thm:Riesz-Gindikin} tells us that
 if a Gindikin-Riesz distribution is a positive measure,
 then
 it equals a Riesz measure associated to the virtual sum 
 $q_\Vsp^{\us} = (q_\Vsp^1)^{\oplus s_1} \oplus \cdots \oplus  (q_\Vsp^r)^{\oplus s_r}$
 of basic quadratic maps up to a constant multiple. 
For example,
 let us recall 
 the \changed{homogeneous $\Pi_r$-positive} quadratic map 
 $q^I : R^I \to \mathrm{Sym}(r,\real)$ with
 \changed{$I \subset \{1, \dots, r\}$}
 and the permutation matrix
 $w_0 \in \mathfrak{S}_r \subset GL(r, \real)$
 \changed{in Section 3.1}.
Putting $g_0 := \rho(w_0)$,
 we have $q^I = g_0 \circ q^{\changed{ \{r-k+1, \dots, r\} } }$
 \changed{($k := \sharp I$)},
 while $q^{\changed{ \{r-k+1, \dots, r\} } }$
 is exactly the basic quadratic map $q^{\changed{ r-k+1 } }_\Vsp$
 for $\Zsp_\Vsp = \mathrm{Sym}(r,\real)$.
Therefore, 
 the Riesz measure $\mu_{(q^I)^{\oplus s}}$ exists if and only if
 $\changed{(0,\dots, 0, \underbrace{s/2, \dots, s/2}_k )} \in \Xi$,
 that is,
 $s \in \{0,1,\dots,k-1\} \cup (k-1,+\infty)$.
In this case,
 $\mu_{(q^I)^{\oplus s}}$ equals the image of 
 $\pi^{s/2} \Rz_{(0,\dots, 0,s/2,\dots,s/2)}
  \changed{= \mu^{}_{(q_\Vsp^{r-k+1})^{\oplus s}}}$
 by $g_0$.   \\[2mm]

Let $q_1 : \real^{m_1} \to \Zsp_\Vsp$
 and $q_2 : \real^{m_2} \to \Zsp_\Vsp$
 be \pg{two} \changed{homogeneous $\Pcone_{\Vsp}$-positive} quadratic maps.
As we have seen in Section 3.1,
 the direct sum $q_1 \oplus q_2$ is not necessarily homogeneous. 
Let us assume that the group 
 $pr_1(\mathrm{Aut}(\Pcone_{\Vsp},q_1)) \cap pr_1(\mathrm{Aut}(\Pcone_{\Vsp},q_2))$
 acts on $\Pcone_\Vsp$ transitively.
In this case,
 we see easily that $q_1 \oplus q_2$ is homogeneous.
As in Proposition~\ref{prop:general_q},
 we can take $g_0 \in G(\Pcone_\Vsp)$ for which
 $g_0 \rho(H_\Vsp)g_0^{-1} \subset pr_1(\mathrm{Aut}(\Pcone_{\Vsp},q_1)) \cap pr_1(\mathrm{Aut}(\Pcone_{\Vsp},q_2))$.
Then we have
 $\det \phi_{q_1}((g_0^{-1})^*\eta) = C_1 \Delta^*_{\um'}(\eta)$
 and
 $\det \phi_{q_2}((g_0^{-1})^*\eta) = C_2 \Delta^*_{\um''}(\eta)$
 for $\eta \in (\real^n)^*$ with some
 $C_1,\,C_2>0$ and $\um',\,\um'' \in \integer^r$.
Now we consider a virtual quadratic map $q = q_1^{\oplus s_1} \oplus q_2^{\oplus s_2}$.
We see that the associated Riesz measure $\mu_q$ exists
 if and only if $s_1 \um' + s_2 \um'' /2 \in \Xi$,
 and in this case, 
 $\mu_q $ is the image of 
 $\changed{C_1^{-s_1/2} C_2^{-s_2/2}}\pi^{(s_1m_1 + s_2 m_2)/2} \Rz_{(s_1 \um'+ s_2 \um'')/2}$ by $g_0$.
Obviously,
 the same argument is valid for general quadratic maps
 $q = q_1^{\oplus s_1} \oplus q_2^{\oplus s_2} \oplus \dots
      \oplus q_t^{\oplus s_t}$.\\    

\noindent{\bf \changed{3.7}. ~Bartlett decomposition of the Wishart laws.}\indent 
Let $q : \real^m \to \Zsp_{\Vsp}$ be a \changed{homogeneous $\Pcone_\Vsp$-positive} quadratic map.
Then the Wishart law $\gamma_{q, \theta}\,\,\,(\theta \in -\Pcone^*_\Vsp)$
 is the image of the normal law 
 $N(0, \phi(-\theta)^{-1})$ on the vector space $\real^m$ by the quadratic map $q/2$,
  \pg{see  Remark \ref{rem-random}}.
 
\pg {However, this description of the Wishart law
  does not permit us to determine  its support   in general}.
In this section,
 we shall give another  \pg{construction} of the Wishart  \pg{random matrices} ,
 which is \changed{a} generalization of the Bartlett decomposition 
 (\changed{\cite{Ba}, \cite[Theorem 3.2.14]{Mu}}) \pg{ and has the advantage of
 controlling the support of the underlying Wishart law.}
Moreover, 
 the result is valid for virtual quadratic maps.

First we consider the virtual quadratic map
 $q_\Vsp^{\us} = (q_\Vsp^1)^{\oplus s_1} \oplus \cdots \oplus  (q_\Vsp^r)^{\oplus s_r}$
\changed{
 whose associated Riesz measure $\mu_{q_\Vsp^{\us}}$ exists.
Then $\usigma = \sum_{i=1}^r s_i \um(i)/2$ belongs to $\Xi$
 and we have 
 $\mu_{q_\Vsp^{\us}} = \pi^{|\usigma|} \Rz_{\usigma}$
 by Theorem~\ref{thm:Riesz-Gindikin}.
Moreover,
 we have
 $L_{\mu_{q_\Vsp^{\us}}}(\changed{\theta}) = \pi^{|\usigma|} \Delta^*_{- \usigma^*}(-\changed{\theta})$.
Therefore we obtain from (\ref{eqn:def_of_Wishart}) that
\begin{equation} \label{eqn:gamma-GR}
  \gamma_{q_\Vsp^{\us}, \theta}(dy) = e^{\cpling{y}\theta}\Delta^*_{\sigma^*}(-\theta) \Rz_{\usigma}(dy) \qquad (y \in \real^n).
\end{equation}
We remark that  distributions of this type are considered in \cite{H-L} for the case  \pg{when}
 $\Pcone_\Vsp$ is a symmetric cone.
}
\changed{
Assume that $\gamma_{q_\Vsp^{\us},\theta}$ is not the Dirac measure.
Then $\usigma \ne (0,\dots,0)$,
 so that
 we can take \pg{a} non-zero $\uep \in \{0,1\}^r$ and $\uu \in R_+(\uep)$
 for which
 $\usigma = \uu + \up(\uep)/2$.
Recall the standard quadratic map 
 $q_\Vsp^{\uep} : W_\Vsp^{\uep} \to \Zsp_\Vsp$
 and the subset $W_{\Vsp}^{\uep, +} \subset W_\Vsp^{\uep}$
 \pg{introduced} in Section 3.4.
As  noted in (\ref{eqn:Tx}),
 each element $x \in W_{\Vsp}^{\uep, +}$ is identified with a lower triangular
 matrix $T_x$ for which
 $q_\Vsp^{\uep}(x) = T_x \,\transp{T_x}$.
Thus,
 the $W_{\Vsp}^{\uep, +}$-valued random variable $X^{\uu}$ in the \pg{following} theorem  
 can be regarded as a triangular random matrix, similarly to the Bartlett decomposition
 of  the classical Wishart distribution.
}


\begin{thm} \label{thm:Bartlett1}
\pg{Let $\usigma = \uu + \up(\uep)/2$  
and $X^{\uu}$}
be an $W_{\Vsp}^{\uep,+}$-valued random variable
 whose components \changed{are independent and }
 satisfy $(X^{\uu}_{ii})^2 \sim \chi^2(2u_i)$ and $X^{\uu}_{li} \sim N(0, I_{n_{li}})$
 for $i \in I(\uep)$ and $l>i$.\\
{\rm (i)} 
The Wishart law $\gamma_{q_\Vsp^{\us}, -I_N}$ is the law of $Y = q_\Vsp^{\uep}(X^{\uu})/2$
\pg{and is supported by  $\orbt_{\uep}$.}\\
{\rm (ii)} 
\changed{For $\theta = - \rho(T)^*I_N \in - \Pcone^*_\Vsp$ with $T \in H_\Vsp$, 
 the Wishart law $\gamma_{q_\Vsp^{\us}, \theta}$} is the law of $Y' = \rho(T)^{-1} \circ q_\Vsp^{\uep}(X^{\uu})/2$
 \pg{and is supported by  $\orbt_{\uep}$.}
\end{thm}
\pf
For a measurable function $f$ on $\Zsp_\Vsp$,
\changed{
 we see from (\ref{eqn:gamma-GR}) and Proposition~\ref{prop:image_meas} that
$$
\int_{\Zsp_\Vsp} f(y) \gamma_{q_\Vsp^{\us}, -I_N}(dy) 
 =  \int_{\Zsp_\Vsp} f(y)e^{-\cpling{y}{I_N}}\Rz_{\usigma}(dy)
 = \int_{W_\Vsp^{\uep,+}} f(q_{\Vsp}^{\uep}(x)) e^{- \cpling{q_{\Vsp}^{\uep}(x)}{I_N} }\,\meas_{\uu}^{\uep}(dx).
$$}
Since $\cpling{q_{\Vsp}^{\uep}(x)}{I_N} = \norm{x}^2$, 
 by the change of variable of $x$ by $x/\sqrt{2}$,
 we rewrite the last term as
$$
 \int_{W_\Vsp^{\uep,+}} f(q_{\Vsp}^{\uep}(x)/2) e^{- \norm{x}^2/2 }\,\meas_{\uu}^{\uep}(dx/2).
$$
Keeping the Remark \ref{rem-aboutM} in mind,
 we see that the law of the random variable $X^{\uu}$ is $e^{- \norm{x}^2/2 }\,\meas_{\uu}^{\uep}(dx/2)$.
Hence (i) holds.
\changed{
To show (ii), 
 it suffices to check that $\gamma_{q_\Vsp^{\us},\theta}$ is the image of 
 $\gamma_{q_\Vsp^{\us},-I_N}$ 
 by $\rho(T)^{-1}$.
Since $\gamma_{q_\Vsp^{\us}, -I_N}(dy) = e^{-\cpling{y}{I_N}} \Rz_{\usigma}(dy)$,
 we have
 \begin{align*}
 \gamma_{q_\Vsp^{\us}, -I_N}(\rho(T)\,dy) = e^{-\cpling{\rho(T)y}{I_N}} \Rz_{\usigma}(\rho(T)\,dy)
 = e^{\cpling{y}{\theta}} \chi_{\usigma}(T)\Rz_{\usigma}(dy)
 \end{align*}
 by (\ref{eqn:GR_T}).
Therefore (\ref{eqn:gamma-GR}) together with 
 (\ref{eqn:def_of_Deltastar})
 leads us to the assertion (ii).
\qed}
$ $\\

\changed{
Now we consider the Wishart distribution $\gamma_{q^{\oplus s}, \theta}$,
 where $q$ is a general homogeneous $\Pcone_\Vsp$-positive quadratic map.
First we show a refinement of the first part of Proposition~\ref{prop:general_q}.}


\begin{lemma} \label{lemma:general_qtheta}
Let $q : \real^m \to \Zsp_\Vsp$ be a homogeneous $\Pcone_\Vsp$-positive quadratic map,
 and $\theta$ an element of  $- \Pcone_{\changed{\Vsp}}^*$.
Then there exist $g_0 \in G(\Pcone_{\changed{\Vsp}})$, 
 $\um \in \integer^r$ and $C >0$
 for which 
\begin{equation} \label{eqn:phi_gq}
 \det \phi_q((g_0^{-1})^*\eta) = C \Delta^*_{\um}(\eta) \qquad (\eta \in \Pcone_\Vsp^{\changed{*}})
\end{equation}
 and $g_0^*\theta = - I_N$.
\end{lemma}

\pf
It is shown in the proof of Proposition~\ref{prop:general_q} that
 there exists $a_0 \in G(\Pcone_\Vsp)$ for which
 $a_0 \rho(H_{\Vsp}) a_0^{-1} \subset pr_1(\mathrm{Aut}(\Pcone_{\Vsp}, q))$.
Since $-a_0^* \theta \in \Pcone^*_\Vsp$,
 we take $T_0 \in H_\Vsp$
 for which
 $\rho^*(T_0)I_N = -a_0^* \theta$.
Put $g_0 := a_0 \rho(T_0)^{-1} \in G(\Pcone_\Vsp)$. 
Then $g_0 \rho(H_{\Vsp}) g_0^{-1} \subset pr_1(\mathrm{Aut}(\Pcone_{\Vsp}, q))$
 and $g_0^* \theta = \rho^*(T_0)^{-1} a_0^* \theta = - I_N$.
Similarly to Proposition~\ref{prop:general_q},
 we see that $g_0$ together with \changed{an} appropriate $\um$ and $C>0$ 
 satisfies the required properties.
\qed


\changed{
Assume that there exists the Riesz measure 
 $\mu_{q^{\oplus s}}$
 associated to a virtual quadratic map $q^{\oplus s}$,
 and that $\mu_{q^{\oplus s}}$ is not the Dirac measure.
Then we can take non-zero 
 $\uep \in \{0,1\}^r$ and $\uu \in R_+(\uep)$
 such that $s\um\changed{/2} = \uu + \up(\uep)/2$
 as in Theorem~\ref{thm:Riesz-Gindikin}. 
Using these data together with $g_0$ in Lemma~\ref{lemma:general_qtheta},
 we obtain the Bartlett decomposition of the Wishart distribution 
 $\gamma_{q^{\oplus s}, \theta}$.
}

\begin{thm} \label{thm:Bartlett2}
\pg{Let  $s\um\changed{/2} = \uu + \up(\uep)/2$ and $X^{\uu}$} be the $W_{\Vsp}^{\uep,+}$-valued random variable in \changed{Theorem~\ref{thm:Bartlett1}}.
Then the Wishart law $\gamma_{q^{\changed{\oplus s}},\theta}$ is the law of $Y = g_0 \circ q_\Vsp^{\uep}(X^{\uu})/2$.   
\end{thm}

\pf
Put $q' :=g_0^{-1} \circ q^{\oplus s}$.
\changed{
As is seen in the proof of Proposition~\ref{prop:general_q},
 the Riesz measure $\mu_{(q')^{\oplus s}}$ equals $C^{-s/2} \pi^{s m/2} \Rz_{s \um/2}$.
Thus,
 similarly to the proof of Theorem~\ref{thm:Bartlett1} (i),
 we see that
 $\gamma_{(q')^{\oplus s}, -I_N}(dy) = e^{-\cpling{y}{I_N}} \Rz_{s \um\pg{/}2}(dy)$,
 and that $\gamma_{(q')^{\oplus s}, -I_N}$ is the law of $Y = q_\Vsp^{\uep}(X^{\uu})/2$.
Since $q^{\oplus \pg{s}} = g_{\pg{0}} \circ (q')^{\oplus s}$,
 Theorem~\ref{thm:Bartlett2} follows from Proposition~\ref{prop:virtual_trans}.
\qed}
$ $\\
 
We have seen that Riesz measures and Wishart laws
 associated to a homogeneous quadratic map
 are obtained
 (up to linear \changed{transforms 
  as in Proposition \ref{prop:general_q} and Theorem~\ref{thm:Bartlett2}}) 
 as the ones associated to a virtual quadratic map
 $q_\Vsp^{\us} = (q_\Vsp^1)^{\oplus s_1} \oplus \cdots \oplus  (q_\Vsp^r)^{\oplus s_r}$,
 that is,
 a virtual sum of basic quadratic maps.
However, 
 it does not mean that every homogeneous quadratic map is  
 equal to a direct sum of basic quadratic maps.
\changed{
The structure of homogeneous quadratic maps is 
 more rich than the maps generated by basic quadratic maps.
}
 \pg{Let us study the following   example.}\\[2mm]
\pg{{\bf Example 5.}} 
Let $\mathrm{Herm}(2, \complex)$ be the vector space of Hermitian matrices of size 2,
 and $\Omega \subset \mathrm{Herm}(2,\complex)$
 the subset of positive definite matrices.
Then we see that
$$
 \Omega = 
 \set{\begin{pmatrix} y_1 & y_3 + \iu y_4 \\ y_3 - \iu y_4 & y_2 \end{pmatrix}}
 {y_1 >0,\,\,y_1 y_2 - (y_3)^2 - (y_4)^2 >0},
$$
 so that $\Omega$ is the $4$-dimensional Lorentz cone.
Recalling (\ref{eqn:Lorentz}),
 we have the linear isomorphism 
\begin{align*}
\iota : \mathrm{Herm}(2, \complex) \owns
\begin{pmatrix} y_1 & y_3 - \iu y_4 \\ y_3 + \iu y_4 & y_2 \end{pmatrix}
\mapsto \begin{pmatrix} y_1 & 0 & y_3 \\ 0 & y_1 & y_4 \\ y_3 & y_4 & y_2 \end{pmatrix} \in \Zsp_\Vsp
\end{align*}
 which gives a matrix realization of $\Omega$.
Let us consider the quadratic map
 $\tilde{q} : \complex^2 \owns z \mapsto z \transp{\bar{z}} \in \mathrm{Herm}(2, \complex)$,
 which is clearly $\Omega$-positive.
We have a group homomorphism
 $$
 \mathrm{GL}(2, \complex) \owns A \mapsto (\tilde{\rho}(A), A) \in \mathrm{Aut}(\Omega, \tilde{q}),
 $$
 where $\tilde{\rho}(A) \in GL(\mathrm{Herm}(2,\complex))$
 is defined by $\tilde{\rho}(A)(Z) := A Z\, \transp{\bar{A}}\,\,\,(Z \in \mathrm{Herm}(2, \complex))$.
Since $\tilde{\rho}(GL(2,\complex))$ acts on $\Omega$ transitively, 
 the quadratic map $\tilde{q}$ is homogeneous.
Keeping the natural isomorphism $\complex^2 \simeq \real^4$ in mind, 
 we define the quadratic map $q: \real^4 \to \Zsp_\Vsp$ by
\begin{align*}
q(x) &:= \iota \circ \tilde{q}\begin{pmatrix} x_1 + \iu x_2 \\ x_3 + \iu x_4 \end{pmatrix}\\
     &= \iota \begin{pmatrix} (x_1)^2 + (x_2)^2 & (x_1 x_3 + x_2 x_4) - \iu (x_1 x_4 - x_2 x_3) \\  
       (x_1 x_3 + x_2 x_4) + \iu (x_1 x_4 - x_2 x_3) & (x_3)^2 + (x_4)^2 \end{pmatrix}\\
     &= \begin{pmatrix} (x_1)^2 + (x_2)^2 & 0 & x_1 x_3 + x_2 x_4 \\ 0 & (x_1)^2 + (x_2)^2 & x_1 x_4 - x_2 x_3\\
          x_1 x_3 + x_2 x_4 & x_1 x_4 - x_2 x_3 & (x_3)^2 + (x_4)^2 \end{pmatrix}.
\end{align*}  
Then we have
$$
\phi_q(\eta) 
 = \begin{pmatrix} \eta_1 & 0 & \eta_3 & \eta_4 \\ 0 & \eta_1 & - \eta_4 & \eta_3 \\
   \eta_3 & - \eta_4 & \eta_2 & 0 \\ \eta_4 & \eta_3 & 0 & \eta_2 \end{pmatrix}
$$
 for $\eta \in \Zsp_\Vsp$.
\changed{
It is easily checked that the map
 $\phi_q \circ \iota : \mathrm{Herm}(2, \complex) \to \mathrm{Sym}(4,\real)$
 is a Jordan algebra representation.
For $\eta \in \Omega^*$ we have
\begin{equation} \label{eqn:q-Laplace}
 L_{\mu_q}(-\eta) = \pi^2 (\det \phi_q(\eta))^{-1/2} = \pi^2 (\eta_1 \eta_2 - (\eta_3)^2 - (\eta_4)^2)^{-1}
\end{equation}
by Lemma \ref{lemma:Laplace_trans}.
}
On the other hand,
 the basic quadratic maps $q_\Vsp^i : W_\Vsp^i \to \Zsp_\Vsp\,\,\,(i=1,2)$ are given by
$$
q_\Vsp^1 \begin{pmatrix} x_1 & 0 \\ 0 & x_1 \\ x_3 & x_4 \end{pmatrix} 
= \begin{pmatrix} (x_1)^2 & 0 & x_1 x_3 \\ 0 & (x_1)^2 & x_1 x_4 \\ x_1 x_3 & x_1 x_4 & (x_3)^2 + (x_4)^2 \end{pmatrix},
 \quad
q_\Vsp^2\begin{pmatrix} 0 \\ 0 \\ x_2 \end{pmatrix} 
= \begin{pmatrix} 0 & 0 & 0 \\ 0 & 0 & 0 \\ 0 & 0 & (x_2)^2 \end{pmatrix},
$$ 
 so that we have for $\eta \in \Zsp_\Vsp^*$ 
$$ 
 \phi_\Vsp^1(\eta) = \begin{pmatrix} \eta_1 & \eta_3 & \eta_4 \\ \eta_3 & \eta_2 & 0 \\ \eta_4 & 0 & \eta_2 \end{pmatrix},
 \quad
 \phi_\Vsp^2(\eta) = \eta_2. 
$$
Thus we obtain
 \begin{equation} \label{eqn:q_12-Laplace}
 L_{\mu_{q_\Vsp^1}}(-\eta) = \pi^{3/2} (\eta_2)^{-1/2}(\eta_1 \eta_2 - (\eta_3)^2 - (\eta_4)^2)^{-1/2}, \quad
 L_{\mu_{q_\Vsp^2}}(-\eta) = \pi^{1/2} (\eta_2)^{-1/2}
 \end{equation}
 for $\eta \in \Pcone_\Vsp^*$.
Comparing (\ref{eqn:q-Laplace}) and (\ref{eqn:q_12-Laplace}),
 we see that
 $\mu_q = \mu_{(q_\Vsp^1)^{\oplus 2} \oplus (q_\Vsp^2)^{\oplus (-2)}}$,
 whereas the quadratic map $q$ is by no means equal to the virtual quadratic map 
 $(q_\Vsp^1)^{\oplus 2} \oplus (q_\Vsp^2)^{\oplus (-2)}$.
We see also from (\ref{eqn:q-Laplace}) and (\ref{eqn:q_12-Laplace}) that
 $\mu_{q  \oplus (q_\Vsp^2)^{\oplus 2}}= \mu_{(q_\Vsp^1)^{\oplus 2}}$,
 whereas
 the two (true) quadratic maps
 $q  \oplus (q_\Vsp^2)^{\oplus 2}$ and $(q_\Vsp^1)^{\oplus 2}$
\pg{ do} not coincide even up to linear transforms $g_0\in G(\Omega)$, \pg{as the domains of
these maps are different}.

Therefore two different quadratic maps may correspond to the same Riesz measure.
\\

%
%
%

\noindent\changed{{\bf 3.8.~Density function for the non-singular case.} \indent
Since the orbit $\orbt_{\uep} = \rho(H)E_{\uep}$
 is contained in the boundary $\partial \Pcone_\Vsp$ 
 of the homogeneous cone $\Pcone_{\Vsp}$ unless $\ep = (1, \dots, 1)$,
 the Gindikin-Riesz distribution $\Rz_{\usigma}$ is a singular measure 
 for $\usigma \in \Xi(\uep)$ with $\uep \ne (1, \dots, 1)$
 thanks to Proposition~\ref{prop:image_meas}.
On the other hand,
 if $\usigma \in \Xi(1, \dots, 1)$,
 that is,
 $$
 \sigma_i > p_i/2 \quad (i=1, \dots, r),
 $$
 where $p_i :=p_i(1, \dots, 1) = \sum_{l>i} n_{li}$,
 then the Gindikin-Riesz distribution is \pg{an} absolutely continuous measure
 with respect to the Lebesgue measure,
 and the density function is given explicitly in \cite{G1} as follows.}

\changed{
Noting that the group $H_\Vsp$ acts on $\Pcone_\Vsp$ simply transitively,
 we define the function $\Delta_{\usigma} : \Pcone \to \complex^{\times}$
 for $\usigma = (\sigma_1, \dots, \sigma_r) \in \complex^r$
 by
 $\Delta_{\usigma}(\rho(T)I_N) := \chi_{\usigma}(T) \,\,\,(T \in H_\Vsp)$.
For $y = \rho(T)I_N = T\, \transp{T} \in \Pcone_\Vsp$,
 we can express $\Delta_{\usigma}(y)$ as a product of powers of principal minors of $y$
 (cf. \cite[p. 122]{F-K}).
Define $\ud = (d_1, \dots, d_r) \in \integer^r/2$ by
 $d_k := 1 + (\sum_{l>k} n_{lk} + \sum_{i<k} n_{ki})/2$.
Then $\Delta_{-\ud}(y)dy$ gives a $G(\Pcone_\Vsp)$-invariant measure on $\Pcone$
 (\cite[Proposition 2.2]{G1}).
Take $\usigma \in \Xi(1, \dots, 1)$.
We see from \cite[Theorem 2.1]{G1} that the integral
 $$ \Gamma_{\Pcone_\Vsp}(\usigma) := \int_{\Pcone_\Vsp} e^{- \cpling{y}{I_N}} \Delta_{\usigma - \ud}(y)\,dy $$
 converges and equals
 $\pi^{(\dim \Zsp_\Vsp - r)/2} \prod_{i=1}^r \Gamma(\sigma_i - p_i/2)$.
By \cite[Proposition 2.3]{G1}, we see that
\begin{equation} \label{eqn:GR_regular}
 \Rz_{\usigma}(dy) = \frac{\Delta_{\usigma - \ud}(y)}{\Gamma_{\Pcone_\Vsp}(\us)}\,dy
 \qquad (y \in \Pcone_\Vsp).
\end{equation}
Owing to (\ref{eqn:gamma-GR}) and (\ref{eqn:GR_regular}),
 we conclude the following proposition.
\begin{prop} \label{prop:density}
Let
 $q_\Vsp^{\us} = (q_\Vsp^1)^{\oplus s_1} \oplus \cdots \oplus  (q_\Vsp^r)^{\oplus s_r}$
 be the virtual quadratic map such that
 $\usigma = \sum_{i=1}^r s_i \um(i)/2$
 belongs to $\Xi(1, \dots, 1)$, 
 that is,
 $\sigma_i > p_i/2 \,\,\,(i=1, \dots, r)$.
Then one has
\begin{equation} \label{eqn:density} 
 \gamma_{q_\Vsp^{\us}, \theta}(dy) 
 = \frac{e^{\cpling{y}\theta}\Delta^*_{\sigma^*}(-\theta)\Delta_{\usigma - \ud}(y)}{\Gamma_{\Pcone_\Vsp}(\us)}\,dy
 \qquad (y \in \Pcone_\Vsp).
\end{equation} 
\end{prop}
}
$ $\\
\pg{ Note that the formula (\ref{eqn:density}) served as 
 a definition of a Wishart law in  \cite{A-W}.}

\noindent \changed{\textbf{Example 6.}
Let $\Zsp_\Vsp$ be the space  \pg{defined} in (\ref{eqn:Vinberg_sp}).
If $y = \rho(T)I_N = T \transp{T} \in \Pcone_\Vsp$ with $T \in H_\Vsp$,
 then we see easily that
$$
 y_{11} = (t_{11})^2,\quad
 y_{11}y_{22} - (y_{21})^2 = (t_{11})^2(t_{22})^2, \quad
 y_{11}y_{33} - (y_{31})^2 = (t_{11})^2(t_{33})^2,
$$ 
 so that $\Delta_{\usigma}(y) = (t_{11})^{2 \sigma_1}(t_{22})^{2 \sigma_2}(t_{33})^{2 \sigma_3}$ equals
$$
 (y_{11})^{\sigma_1 - \sigma_2 - \sigma_3}(y_{11}y_{22} - (y_{21})^2)^{\sigma_2} (y_{11}y_{33} - (y_{31})^2)^{\sigma_3}.
$$
On the other hand,
 we have $(p_1,\, p_2,\, p_3) = (0,1,1)$ and $(d_1,\,d_2,\,d_3) = (2,\,3/2,\,3/2)$.
Thus we have by (\ref{eqn:GR_regular})
 $$ \Rz_{(\sigma_1,\sigma_2,\sigma_3)}(dy) 
  = \frac{(y_{11})^{\sigma_1 - \sigma_2 - \sigma_3+1}(y_{11}y_{22} - (y_{21})^2)^{\sigma_2-3/2} (y_{11}y_{33} - (y_{31})^2)^{\sigma_3-3/2}}
         {\pi \Gamma(\sigma_1)\Gamma(\sigma_2 - 1/2) \Gamma(\sigma_3 - 1/2)}\,dy $$
 if $\sigma_1 > 0,\,\sigma_2>1/2$ and $\sigma_3 >1/2$.
Let us consider the Wishart laws associated to
 the virtual quadratic map $(q_{\Vsp}^1)^{\oplus s}$,
 where $q_\Vsp^1: W_\Vsp^1 \to \Zsp_\Vsp$ is the basic quadratic map.
Since $\um(1) = (1,\,1,\,1)$,
 we observe that $s \um(1)/2 \in \Xi$ if and only if
 $s \in \{0,1\} \cup (1,+\infty)$.
If $s=0$, the associated Wishart law is the Dirac measure.
If $s=1$, the associated Wishart law $\gamma_{q_\Vsp^1,\,\theta}\,\,\,(\theta  \in -\Pcone_\Vsp^*)$ is described
 as the image of the normal law $N(0, \phi_\Vsp^1(-\theta)^{-1})$ on $W_\Vsp^1 \equiv \real^3$ by the quadratic map
 $q_\Vsp^1/2$,
 where $\phi_\Vsp^1(-\theta)$ is given in Example 4 after Proposition~\ref{prop:dualcone}.
}

\changed{
If $s>1$, then $s \um(1)/2 = (s/2,\,s/2,\,s/2)$ belongs to $\Xi(1,1,1)$.
Since 
 $$ \Delta^*_{s \um(1)^*/2}(\eta) = \det \phi_{\Vsp}^1(\eta)^{s/2} = \bigl( \eta_{11} \eta_{22} \eta_{33} - \eta_{33}(\eta_{21})^2 - \eta_{22}(\eta_{31})^2 \bigr)^{s/2}
   \qquad (\eta \in \Pcone_\Vsp^*),$$
 we have for $\theta = - \eta \in - \Pcone_\Vsp^*$  
 \begin{align*}
 \gamma_{(q_\Vsp^1)^{\oplus s},\,\theta}(dy) 
 &= \frac{e^{-\cpling{y}{\eta}} \bigl( \eta_{11} \eta_{22} \eta_{33} - \eta_{33}(\eta_{21})^2 - \eta_{22}(\eta_{31})^2 \bigr)^{s/2}}
        {\pi \Gamma(s/2)\Gamma((s - 1)/2) \Gamma((s- 1)/2)} \times\\
 & \quad \times (y_{11})^{1-s/2}(y_{11}y_{22} - (y_{21})^2)^{(s-3)/2} (y_{11}y_{33} - (y_{31})^2)^{(s-3)/2} \,dy 
 \quad (y \in \Pcone_\Vsp)
\end{align*}
 by Proposition~\ref{prop:density}.
}

\end{document}